\documentclass{article}
\usepackage{amsmath,amssymb}
\usepackage[all]{xy}
\newtheorem{thm}{Theorem}[section]
\newtheorem{prop}[thm]{Proposition}
\newtheorem{cor}[thm]{Corollary}
\newtheorem{lem}[thm]{Lemma}
\newenvironment{dfn}{\medskip\refstepcounter{thm}
\noindent{\bf Definition \thesection.\arabic{thm}\ }}{\medskip}

\newenvironment{proof}[1][,]{\medskip\ifcat,#1
\noindent{{\it Proof}:\ }\else\noindent{\it Proof of #1.\ }\fi}
{\hfill$\square$\medskip}
\newenvironment{note}[1][Note]{\begin{trivlist}
\item[\hskip \labelsep {\bfseries #1}]}{\end{trivlist}}
\newenvironment{notes}[1][Notes]{\begin{trivlist}
\item[\hskip \labelsep {\bfseries #1}]}{\end{trivlist}}
\newenvironment{remark}[1][Remark]{\begin{trivlist}
\item[\hskip \labelsep {\bfseries #1}]}{\end{trivlist}}
\newenvironment{remarks}[1][Remarks]{\begin{trivlist}
\item[\hskip \labelsep {\bfseries #1}]}{\end{trivlist}}
\newenvironment{ack}[1][Acknowledgements]{\begin{trivlist}
\item[\hskip \labelsep {\bfseries #1}]}{\end{trivlist}}
\newenvironment{funding}[1][Funding]{\begin{trivlist}
\item[\hskip \labelsep {\bfseries #1}]}{\end{trivlist}}
\def\eq#1{{\rm(\ref{#1})}}
\def\H{{\mathbb H}}
\def\R{{\mathbb R}}
\def\Z{{\mathbb Z}}
\def\N{{\mathbb N}}
\def\O{{\mathbb O}}
\def\SS{{\mathbb S}}
\DeclareMathOperator\GL{GL}

\DeclareMathOperator\SO{SO}
\DeclareMathOperator\U{U}
\DeclareMathOperator\Sp{Sp}
\DeclareMathOperator\id{id}
\def\Im{\mathop{\rm Im}\nolimits}
\def\Re{\mathop{\rm Re}\nolimits}
\def\d{{\rm d}}
\def\w{\wedge}
\def\C{{\mathbb C}}
\DeclareMathOperator\Spin{Spin}
\DeclareMathOperator\GG{G}
\DeclareMathOperator\ind{ind}
\DeclareMathOperator\Ker{Ker}

\DeclareMathOperator\vol{vol}
\def\G2{\GG_2}
\newcommand{\Dslash}{\ensuremath \hspace{1pt}\raisebox{0.025cm}{\slash}\hspace{-0.24cm} D}
\begin{document}

\title{Asymptotically Conical Associative 3-folds}
\author{\textsc{Jason D. Lotay}\\ Mathematical Sciences Research Institute, Berkeley}
\date{}
\maketitle \linespread{1.1}

\begin{center}
{\large\bf Abstract}
\end{center}

\noindent Given an associative 3-fold $N$ in $\R^7$ which is asymptotically conical with generic rate $\lambda<1$, 
we show that the moduli space of deformations of $N$ is locally 
homeomorphic to the kernel of a smooth map between smooth manifolds.  Moreover, the virtual dimension of the moduli 
space is computed and shown to be non-negative for $\lambda>-1$, whereas $N$ is expected to be isolated for $\lambda\leq -1$.

\section[Introduction]{Introduction}\label{s1}

Harvey and Lawson defined a distinguished class of 3-dimensional submanifolds of $\R^7$, known as 
\emph{associative} 3-folds, in their seminal paper on calibrated geometries \cite{HarLaw}.  
The definition of these submanifolds can be found in $\S$\ref{s2}.  They can be thought of as generalisations 
of \emph{special Lagrangian} 3-folds in $\C^3$ which have been well-studied by a number of authors.  
Explicit examples of associative 3-folds which do not come from lower-dimensional geometries
were constructed by the author in \cite{Lotay3} and \cite{Lotay2}.  
In particular, a family was produced consisting of \emph{asymptotically conical} (AC) associative 3-folds; 
that is, associative 3-folds which converge to a cone at infinity with a given 
decay rate.  This family is presented in $\S$\ref{s7} and is part of the principal inspiration for this paper.  

The other main motivations come from the deformation theory studies of McLean \cite[$\S$5]{McLean} on \emph{compact} 
associative 3-folds, and Marshall \cite{Marshall} and the author \cite{Lotay1} on AC special Lagrangian and
\emph{coassociative}
submanifolds respectively.  Coassociative submanifolds are examples of calibrated 4-dimensional 
submanifolds of $\R^7$ and were also defined in \cite{HarLaw}.   
Asymptotically conical special Lagrangian submanifolds are also discussed in \cite{Pacini} and in the series of papers 
\cite{Joyce5}-\cite{Joyce9}.  Asymptotically cylindrical deformations of calibrated submanifolds
 are studied in  
\cite{JoySalur}, \cite{Marshall} and \cite{Salur}.  Recently, Gayet and Witt \cite{Gayet} 
described the deformation theory of associative 3-folds with boundary in a coassociative 4-fold.

From McLean's study \cite[$\S$5]{McLean}, 
we cannot say, in general, whether there is a smooth moduli space of compact associative deformations and, even if there is one, 
its virtual dimension is zero.  Thus, we would generically expect compact associative 3-folds to be isolated. 
For the AC setting, we still 
cannot determine generally whether the moduli space of deformations is smooth, though we give a condition for 
a smooth moduli space to exist.  
However, in contrast to the compact case, we show that the expected dimension of the
 moduli space is guaranteed to be non-negative for certain decay rates and is even positive if the 
link of the asymptotic cone has certain properties -- a rather surprising result.

\medskip

We prove two main theorems.  The first is the following.

\begin{thm}\label{thm1}
Let $N$ be an associative 3-fold in $\R^7$ which is asymptotically conical with generic rate $\lambda<1$ to a cone $C$ 
in $\R^7$ with link $\Sigma$ in $\mathcal{S}^6$, as defined in Definition \ref{s2dfn8}. 

The moduli space $\mathcal{M}(N,\lambda)$, given in Definition \ref{s4dfn1}, 
of associative deformations of $N$ which are asymptotically conical 
to $C$ with rate $\lambda$ is locally homeomorphic to the 
kernel of a smooth map $\pi$ between smooth manifolds.  The manifolds are open neighbourhoods of zero in
finite-dimensional Banach spaces and $\pi(0)=0$.  Therefore, $\mathcal{M}(N,\lambda)$ is a smooth manifold if\/ $\pi$
is the zero map.
\end{thm} 

Our second main theorem gives the expected dimension of the moduli space.

\begin{thm}\label{thm2}  Use the notation of Theorem \ref{thm1}. 
 There exists an elliptic first-order differential operator $\bar{\partial}_\Sigma$ on the normal bundle of\/ $\Sigma$ in $\mathcal{S}^6$ such that, if 
    $$\d(\mu)=\dim\Ker\big(\bar{\partial}_\Sigma-(\mu+1)\id\!\big)$$ for $\mu\in\R$,
then the expected dimension of $\mathcal{M}(N,\lambda)$ for $\lambda>-1$ is
\begin{equation}\label{moddimeq1}
\frac{1}{2}\,\d(-1)+\!\!\!\!\sum_{\mu\in(-1,\lambda)}\!\!\!\!\d(\mu).
\end{equation}

Furthermore, the virtual dimension of $\mathcal{M}(N,\lambda)$ is less than or equal to
$-\frac{1}{2}\,\d(-1)$ for $\lambda\leq -1$, and is strictly negative for $\lambda<-2$.
\end{thm}

\begin{remarks}
\begin{itemize}\item[]
\item[(a)]
By a `generic rate $\lambda$' in Theorem \ref{thm1} we mean that $\lambda+1$ does not lie in the  
 spectrum of the operator $\bar{\partial}_\Sigma$ 
mentioned in Theorem \ref{thm2}.  The word `generic' is justified because this spectrum is a countable 
discrete subset of $\R$ as $\bar{\partial}_\Sigma$ is elliptic and $\Sigma$ is compact.  This fact also ensures that the sum in \eq{moddimeq1} is well-defined.
\item[(b)] Any translation of an AC associative 3-fold in $N$ in $\R^7$ is trivially a deformation of $N$ which is still AC for 
rates $\lambda\in(0,1)$, so we would naturally expect the moduli space to have positive dimension for growth rates in this range.
\item[(c)] We show that $\frac{1}{2}\,\d(-1)$ is an integer in Lemma \ref{s6lem5}.
\end{itemize}
\end{remarks}

\medskip

We begin in $\S$\ref{s2} with the basic definitions we require.  We also describe an isomorphism between the normal
bundle of an associative 3-fold and a twisted spinor bundle.  This is followed, in $\S$\ref{s3}, with a 
summary of the relevant theory of weighted Banach spaces, which provides the analytic framework for the paper.   

In $\S$\ref{s4} we set up our deformation problem by characterising local deformations of AC associative 3-folds as  
elements of the kernel of a map $F$ on the normal bundle. 
This map is a nonlinear first-order differential operator whose linearisation at zero, denoted by $\Dslash$, 
 can be thought of as a Dirac operator on twisted spinors.
We also discuss the regularity of elements of $F^{-1}(0)$.  
The next section, $\S$\ref{s5}, is dedicated to describing the 
Fredholm and index theory of $\Dslash$.

The key results leading to Theorems \ref{thm1} and \ref{thm2} are presented in $\S$\ref{s6}.  We identify the infinitesimal 
deformation and obstruction spaces, and  
prove that the moduli space, near $N$, is homeomorphic 
to the kernel, near zero, of a smooth map from an open subset in the
 infinitesimal deformation space to the obstruction space.  We also calculate the expected dimension of the moduli
 space.

Finally, in $\S$\ref{s7}, we discuss examples of AC associative 3-folds.

\begin{notes}\begin{itemize}\item[]
\item[(a)] Manifolds are taken to be nonsingular and submanifolds to
be embedded, for convenience, unless stated otherwise.
\item[(b)] We use the convention that the natural numbers
$\N=\{0,1,2,\ldots\}$.
\item[(c)] Given a manifold $M$ and a vector bundle $E$ over $M$, we write
$\Omega^m(M;E)$ for the bundle of $m$-forms over $M$ taking values in $E$.
\end{itemize} \end{notes}

\section[Asymptotically conical associative submanifolds of R7]{Asymptotically conical associative\\ submanifolds of {\boldmath $\R^7$}}\label{s2}

Associative 3-folds form part of the \emph{calibrated geometry} in seven dimensions associated with the Lie group $\G2$.
We therefore begin by defining calibrations and calibrated submanifolds following the approach in 
\cite[Definition 4.1.3]{Joyce3}.  

\begin{dfn}\label{s2dfn1}
Let $(M,g)$ be a Riemannian manifold.  An \emph{oriented tangent\
$k$-plane} $V$ on $M$ is an oriented $k$-dimensional vector subspace
$V$ of $T_xM$, for some $x$ in $M$. Given an oriented tangent
$k$-plane $V$ on $M$, $g|_V$ is a Euclidean metric on $V$
and hence, using $g|_V$ and the orientation on $V$, we have a
natural volume
form, $\text{vol}_V$, which is a $k$-form on $V$. 

Let $\eta$ be a closed $k$-form on $M$.  Then $\eta$ is a
\emph{calibration} on $M$ if $\eta | _V = \kappa
\cdot\text{vol}_V$ for some $\kappa\leq 1$, for all
oriented tangent $k$-planes $V$ on $M$.

Let $N$ be an oriented $k$-dimensional submanifold of $M$.  Then $N$
is a \emph{calibrated submanifold} or \emph{$\eta
$-submanifold} if $\eta | _{T_xN} =\,\text{vol}_{T_xN}$ for all $x\in N$.
\end{dfn}

Calibrated submanifolds are examples of \textit{minimal} submanifolds \cite[Theorem II.4.2]{HarLaw}.  
  We now define \emph{associative 3-folds} by introducing a distinguished 3-form on $\R^7$ as in 
\cite[Definition 12.1.1]{Joyce3}.

\begin{dfn}\label{s2dfn2}
Let $(x_1,\ldots,x_7)$ be coordinates on $\R^7$.  We
shall denote $\d x_i\w\d
x_j\w\d x_k$ by $\d{\bf x}_{ijk}$.  Define a 3-form $\varphi$ on $\R^7$ by
\begin{equation}\label{s2eq1}
\varphi = \d{\bf x}_{123}+\d{\bf x}_{145}
+\d{\bf x}_{167}+\d{\bf x}_{246}-\d{\bf x}_{257}
-\d{\bf x}_{347}-\d{\bf x}_{356}.
\end{equation}
By \cite[Theorem IV.1.4]{HarLaw}, $\varphi$ is a calibration on $\R^7$ and 
$\varphi$-submanifolds are known as \emph{associative 3-folds}.  The stabilizer of $\varphi$ in $\GL(7,\R)$ is $\G2$.
\end{dfn}

\begin{note}
Since $*\varphi$ is also a calibration \cite[Theorem IV.1.16]{HarLaw}, we can define $*\varphi$-submanifolds of $\R^7$:
 these are called \emph{coassociative 4-folds}.
\end{note}

Other well known examples of calibrated submanifolds are given by complex submanifolds of K\"ahler manifolds, and 
\emph{special Lagrangian submanifolds} of Calabi--Yau manifolds.  We define the latter in the special case when 
the Calabi--Yau manifold is $\C^3$ as they shall be discussed in $\S$\ref{s7}.  

\begin{dfn}\label{SLdfn}
Let $\C^3$ 
be endowed with the standard metric, K\"ahler form $\omega$ and 
holomorphic volume form $\Omega$.  A (real) 3-dimensional submanifold $L$ of $\C^3$ is a
\emph{special Lagrangian 3-fold} if $L$ is calibrated with respect to $\text{Re}\,\Omega$.  
\end{dfn}

As noted in \cite[$\S$12.2.1]{Joyce3}, associative 3-folds generalise calibrated geometries in $\C^3$ in the 
following sense.

\begin{lem}\label{s2lem1}
Decompose $\R^7$ as $\R\oplus\C^3$.  
\begin{itemize}
\item[\emph{(i)}] $\Sigma$ is a holomorphic curve in $\C^3$ if and only if\/ $\R\times\Sigma$ is associative in $\R^7$.
\item[\emph{(ii)}] $L$ is a special Lagrangian 3-fold in $\C^3$ if and only if\/ $\{x\}\times L$ is associative 
in $\R^7$ for $x\in\R$.
\end{itemize}
\end{lem}

\noindent These observations follow easily from the formula $\varphi=dx_1\w\omega+\text{Re}\,\Omega$ on $\R^7=\R\oplus\C^3$, 
where $x_1$ is the coordinate on the $\R$ factor.

\medskip

We will also need to understand some aspects of the geometry of $\mathcal{S}^6$, so we make the following definition.

\begin{dfn}\label{6spheredfn}
Let $\Im\O$ denote the imaginary octonions and identify the standard orthonormal basis of $\R^7$ with 
a basis $\{e_1,\ldots,e_7\}$ for $\Im\O$ such that the cross product on $\Im\O$ 
is given by:
\begin{equation}\label{crossprodeq}
\varphi(x,y,z)=g_{\R^7}(x\times y,z)\qquad\text{for $x,y,z\in\Im\O$,}
\end{equation}
where $\varphi$ is given in \eq{s2eq1} and $g_{\R^7}$ is the Euclidean metric on $\R^7$.

We therefore have a natural embedding $\mathcal{S}^6\hookrightarrow\Im\O$ and we may
identify $T_x\mathcal{S}^6$ with the (real) 6-plane in $\Im\O$ orthogonal to $x\in\mathcal{S}^6$.  Thus, for all $x\in\mathcal{S}^6$, 
we may define a map $J_x:T_x\mathcal{S}^6\rightarrow T_x\mathcal{S}^6$ by $J_x(u)=x\times u$.  
Then $J$ defines an \emph{almost complex structure} on $\mathcal{S}^6$ which is \emph{not} integrable.  
Moreover, if $g_{\mathcal{S}^6}$ is the round metric on $\mathcal{S}^6$, we may define a 2-form $\omega$ on $\mathcal{S}^6$ by $\omega(u,v)=g_{\mathcal{S}^6}(Ju,v)$.  This 2-form is nondegenerate but it is \emph{not} closed.
\end{dfn}

We shall find it useful to have an alternative characterisation of associative 3-folds, given
by \cite[Corollary IV.1.7]{HarLaw}, for which we must first define a vector-valued 3-form on $\R^7$.

\begin{dfn}\label{s2dfn3}
Define $\chi\in C^{\infty}\big(\Omega^3(\R^7;T\R^7)\big)$ by:
\begin{equation}\label{s2eq2}
g_{\R^7}\big(\chi(x,y,z),w\big)=*\varphi(x,y,z,w),
\end{equation}
where $g_{\R^7}$ is the Euclidean metric on $\R^7$, $\varphi$ is given in \eq{s2eq1}, $*$ is the Hodge star with respect to $g_{\R^7}$ and
$x$, $y$, $z$, $w \in C^{\infty}(T\R^7)$.  We can write $\chi$ explicitly in coordinates on $\R^7$:
\begin{align}
\chi&=\nonumber\\
&(-\d\mathbf{x}_{357}+\d\mathbf{x}_{346}+\d\mathbf{x}_{256}+\d\mathbf{x}_{247})\mathbf{e}_1 
+(-\d\mathbf{x}_{367}-\d\mathbf{x}_{345}-\d\mathbf{x}_{156}-\d\mathbf{x}_{147})\mathbf{e}_2 
\nonumber\\
&+(\d\mathbf{x}_{267}+\d\mathbf{x}_{245}+\d\mathbf{x}_{157}-\d\mathbf{x}_{146})\mathbf{e}_3 
+(\d\mathbf{x}_{567}-\d\mathbf{x}_{235}+\d\mathbf{x}_{136}+\d\mathbf{x}_{127})\mathbf{e}_4 
\nonumber\\
&+(\d\mathbf{x}_{467}+\d\mathbf{x}_{234}-\d\mathbf{x}_{137}+\d\mathbf{x}_{126})\mathbf{e}_5 
+(-\d\mathbf{x}_{457}-\d\mathbf{x}_{237}-\d\mathbf{x}_{134}-\d\mathbf{x}_{125})\mathbf{e}_6 
\nonumber\\
&+(\d\mathbf{x}_{456}+\d\mathbf{x}_{236}+\d\mathbf{x}_{135}-\d\mathbf{x}_{124})\mathbf{e}_7, \label{s2eq3}
\end{align}
where $\mathbf{e}_j=\frac{\partial}{\partial x_j}$.
\end{dfn}

\begin{prop}\label{s2prop4}
A 3-dimensional submanifold $N$ of\/ $\R^7$ is associative if and only if\/ 
$\chi|_N\equiv 0$ and $\varphi|_N>0$.
\end{prop}

\begin{remarks}
\begin{itemize}\item[]
\item[(a)] The condition $\chi|_N\equiv 0$ forces $\varphi$ to be a non-vanishing 3-form on $N$, so the positivity of $\varphi$ amounts
to a choice of orientation on $N$.
\item[(b)] The reader may notice that the system of equations given by $\chi|_N\equiv 0$ appears to be overdetermined.  
Further, it can be shown that one may equivalently define associative 3-folds by the vanishing of $\varphi$ on the 
normal bundle of $N$ in $\R^7$, which now is a determined system.  Therefore it would be perhaps more natural to use
 this alternative definition. 
However, we make our choice so as to utilise the framework provided by McLean \cite{McLean} 
in his discussion of deformations of associative 3-folds.
\end{itemize}
\end{remarks}

We need some elementary properties of $\chi$, which follow from explicit 
calculations involving the formula \eq{s2eq3}.

\begin{prop}\label{s2prop5} If $N$ is an oriented 3-dimensional submanifold of $\R^7$,
$\chi|_N\in C^{\infty}\big(\Omega^3(N;\nu(N))\big)$.  If $N$ is associative, $v\mapsto (v\cdot\chi)|_{N}$
 defines an injective map from
the normal bundle $\nu(N)$ of $N$ into $\Omega^2(N;\nu(N))$.
\end{prop}

We shall need to differentiate normal vector-valued forms later on, so we make the following convenient definition.

\begin{dfn}\label{normaldiffdfn}
Let $P$ be an oriented submanifold of a Riemannian manifold $M$.  The Levi--Civita connection on $M$ induces a connection
 $\nabla^{\perp}_P$ on the normal bundle of $P$ in $M$ which takes sections of $\nu(P)$ to sections of $T^*P\otimes\nu(P)$.  Using the 
 connection $\nabla^{\perp}_P$, we may define a derivative $d^{\perp}_P$ acting on $\nu(P)$-valued $m$-forms $\eta$ on $P$ so that $d^{\perp}_P\eta$ is a 
 $\nu(P)$-valued $(m+1)$-form on $P$.  
\end{dfn}

We now describe how to think of normal vector fields as twisted spinors, as in \cite[p.~734-735]{McLean}.

\begin{dfn}\label{s2dfn6}
Let $N$ be an associative 3-fold in $\R^7$.  Since $N$ is spin, we can choose a spin structure on $N$ and let
$\SS(N)$ denote the spinor bundle.  Further, since $N$ is oriented, we can use the orientations on $N$ and $\R^7$ to define 
an orientation on $\nu(N)$.  We let $\mathcal{F}(N)$ denote the bundle of oriented orthonormal bases of the fibres of 
$\nu(N)$.  Then $\mathcal{F}(N)$ is a principal $\SO(4)$ bundle over $N$.  
Moreover, the spin structure on $N$ 
 determines a principal $\Spin(4)=\Sp(1)\times\Sp(1)$ bundle $\mathcal{P}(N)$ which is 
a double cover of $\mathcal{F}(N)$. 

Let $\mathcal{E}(N)$ be the vector bundle over $N$ associated to $\mathcal{P}(N)$ 
via the representation $\epsilon:\Sp(1)\times\Sp(1)
 \rightarrow \GL(\H)$ given by $\epsilon(p,q)(r)=pr$ for $r\in\H$, where $\H$ denotes the quaternions.
 We notice that $\SS(N)$ and $\nu(N)$ are associated to 
$\mathcal{P}(N)$ via the representations $\sigma(p,q)(r)=r\bar{q}$ and
 $\tau(p,q)(r)=pr\bar{q}$ respectively.  Thus, there is an isomorphism 
$\jmath_N:\nu(N)\rightarrow \mathcal{T}(N):=\SS(N)\otimes_{\H}\mathcal{E}(N)$; that is,
the normal bundle is isomorphic to a twisted spinor bundle.
\end{dfn} 

\begin{note}  We can locally interpret the isomorphism $\jmath_N$ as 
follows.
 Recall the notation introduced in Definition \ref{6spheredfn}.  At $x\in N$ identify $T_x\R^7$ with $\Im\O$ such that $T_xN=\langle e_1,e_2,e_3\rangle$.  This can be done smoothly 
since $\G2$ acts transitively on the set of associative 3-planes by \cite[Proposition 12.1.2]{Joyce3}. 
 Therefore, we may locally
view a section $v$ of $\nu(N)$ as a twisted spinor by writing $\jmath_N(v)=v_4e_4+v_5e_5+v_6e_6+v_7e_7$.
\end{note}

Before discussing AC submanifolds we clarify what we mean by a cone and discuss the relationship of 
associative cones in $\R^7$ with the geometry of $\mathcal{S}^6$.

\begin{dfn}\label{s2dfn7}
A \emph{cone} $C$ in $\R^n$ is a submanifold which is nonsingular, except
possibly at 0, and satisfies $e^tC=C$ for all $t\in\R$.  We say that $\Sigma=C\cap\mathcal{S}^{n-1}$
is the \emph{link} of $C$.  
\end{dfn}

We now make an elementary observation.

\begin{lem}\label{phololem}
Use the notation of Definition \ref{6spheredfn}.  A cone $C$ in $\R^7$ is associative if and only if 
its link $\Sigma$ in $\mathcal{S}^6$ is a \emph{pseudoholomorphic curve}; that is, $J(T_\sigma\Sigma)=T_\sigma\Sigma$ 
for all $\sigma\in\Sigma$ or, equivalently, $\omega|_\Sigma=\vol_\Sigma$.
\end{lem}
\begin{proof}
Consider $\R^7\setminus\{0\}\cong\R^+\times\mathcal{S}^6$, with $r$ the 
radial coordinate.  Then, we may write
\begin{align}\label{phiconeeq}
\varphi=\frac{1}{3}\,r^3d\omega+r^2 dr\wedge\omega,
\end{align}
where $\varphi$ is given in \eq{s2eq1}.  The result follows from the definitions.
\end{proof}

\begin{note}
Using \eq{phiconeeq}, we deduce that
\begin{equation*}
*_{\R^7}\varphi=\frac{1}{2}\,r^4\omega\wedge\omega-\frac{1}{3}\,r^3 dr\wedge*_{\mathcal{S}^6}d\omega,
\end{equation*}
where $*_{\R^7}$ and $*_{\mathcal{S}^6}$ are the Hodge stars on $\R^7$ and 
$\mathcal{S}^6$ respectively.  Thus, there exist a 3-form $\xi$ and 2-form $\zeta$ on 
$\mathcal{S}^6$ taking values in $T\R^7|_{\mathcal{S}^6}$ such that
\begin{equation}\label{chiconeeq}
\chi=r^2\xi+rdr\wedge\zeta.
\end{equation}
\end{note}

\begin{dfn}\label{s2dfn8}
Let $N$ be a closed 3-dimensional
submanifold of $\R^7$. Then $N$ is \emph{asymptotically conical}
(AC), or an \emph{AC 3-fold}, \emph{with rate $\lambda$} if there
exist constants $\lambda<1$ and $R>1$, a compact subset $K$ of $N$,
a compact 2-dimensional Riemannian submanifold $(\Sigma,h)$ of
$\mathcal{S}^6\subseteq\R^7$, where $h$ is the restriction of the
round metric on $\mathcal{S}^6$ to $\Sigma$, and a diffeomorphism
$\Psi:(R,\infty)\times\Sigma\rightarrow N\setminus K$ satisfying 
\begin{equation}\label{s2eq4}
\Psi(r,\sigma)-\iota(r,\sigma)\in \big(T_{r\sigma}C\big)^{\perp}\qquad\text{for all 
$(r,\sigma)\in (R,\infty)\times\Sigma$}\end{equation}
and
\begin{equation}\label{s2eq5}
\big|\nabla^j\big(\Psi(r,\sigma)-\iota(r,\sigma)\big)\big|=O\big(r^{\lambda-j}\big)
\qquad\text{for $j\in\N$ as $r\rightarrow\infty$,}
\end{equation}
\noindent where $\iota(r,\sigma)=r\sigma$, $\nabla$ is the
Levi--Civita connection of the cone metric $g_\text{cone}=dr^2+r^2h$
on $C=\iota((0,\infty)\times\Sigma)$ coupled with partial differentiation
on $\R^7$, and
 $|\,.\,|$ is calculated with respect to $g_\text{cone}$.

We say
that $N$ is asymptotically conical to $C$ for clarity when
necessary and call the components of $N_{\infty}=N\setminus K$ the \emph{ends} of $N$.
\end{dfn}

\begin{remarks} 
\begin{itemize}\item[]
\item[(a)]
 The condition $\lambda<1$ in the definition above ensures that, by \eq{s2eq5}, 
the restriction of the Euclidean metric on $\R^7$
to $N$ converges to $g_{\text{cone}}$ at infinity.  
\item[(b)] The condition \eq{s2eq4} simply ensures that $\Psi$ is determined uniquely.
\end{itemize}
\end{remarks}

\begin{dfn}\label{s2dfn9}
Let $N$ be an AC 3-fold in $\R^7$ and use the notation of
Definition \ref{s2dfn8}.  A \emph{radius function}
$\rho:N\rightarrow [1,\infty)$ on $N$ is a smooth map satisfying
$c_1r<\Psi^*(\rho)<c_2r$ for some positive constants $c_1<1$ and $c_2>1$.
\end{dfn}

\noindent We may define a radius function
$\rho$ on $N$ by setting $\rho =1$ on $K$,
$\rho\big(\Psi(r,\sigma)\big)=r$ for $r>R+1$, and then
extending $\rho$ smoothly to our required function on $N$.

\medskip

We conclude the section with an elementary result which is analogous to
\cite[Proposition 2.8]{Lotay1}, and can be proved in a similar manner using Proposition \ref{s2prop4}.  

\begin{prop}\label{s2prop10}
Suppose that $N$ is an associative 3-fold in $\R^7$ which is AC with rate $\lambda$ 
to a cone $C$ in $\R^7$.  Then $C$ is associative.
\end{prop}

\begin{note}
Suppose $N$ is an associative 3-fold in $\R^7$ which is AC with rate $\lambda$ to $C$.  Notice that 
the dilation $e^tN$, for $t\in\R$, is also an associative 3-fold which is AC with rate $\lambda$ to $C$.  Therefore 
$N$ always has a one-parameter family of deformations given by dilation, unless $N$ is a union of 3-planes.   
\end{note}

\section[Weighted Banach spaces]{Weighted Banach spaces}\label{s3}

Here we briefly review the theory of weighted Banach spaces, following the definitions in \cite[$\S$1]{Bartnik}. 
For this section, let $N$ be an asymptotically conical 3-fold in $\R^7$, let $g$ be the metric on $N$ and let 
$\rho$ be a radius function on $N$.  We also let $E$ be a vector bundle on $N$ which is equipped with 
Euclidean metrics on its
fibres and endowed with a connection $\nabla$ that preserves these metrics.
 
We denote Sobolev and H\"older spaces of sections of $E$ by $L^p_k(E)$ and $C^{k,\,a}(E)$, for $p\geq 1$, $k\in\N$ and $a\in(0,1)$.  
The definitions of these Banach spaces can be found in \cite[$\S$1.2]{Joyce3}.  We also 
add the `loc' subscript to our notation to indicate that the sections `locally' lie in the Banach space; 
for example, $C^k_{\text{loc}}(E)$ is the set of sections $\xi$ of $E$ such that $f\xi\in C^k(E)$
for all smooth compactly supported functions $f$ on $N$.

\begin{dfn}\label{s3dfn1}
Let $p\geq 1$, $k\in\N$ and $\mu\in\R$.  The \emph{weighted Sobolev space}
$L_{k,\,\mu}^p(E)$ is the 
subspace of $L_{k,\,\text{loc}}^p(E)$ such that the norm
\begin{equation*}
\|\xi\|_{L_{k,\,\mu}^p}=\left(\sum_{j=0}^k\int_N
|\rho^{j-\mu}\nabla^j\xi|^p\rho^{-3} \,dV_g\right)^\frac{1}{p}
\end{equation*}
is finite.  These are Banach spaces and are also Hilbert spaces when $p=2$.  
Note that $L^p=L^p_{0,\,-\frac{3}{p}}$.
\end{dfn}

\begin{dfn}\label{s3dfn2}
Let $\mu\in\R$ and let $k\in\N$.  The
\emph{weighted $C^k$-space} $C_{\mu}^{k}(E)$ 
is the subspace of $C^k_{\text{loc}}(E)$ such that the norm 
$$\|\xi\|_{C_{\mu}^{k}}=\sum_{j=0}^k \sup_N|\rho^{j-\mu}\nabla^j\xi|$$
is finite.  We also define $C_{\mu}^{\infty}(E)=\cap_{k\geq
0}C_{\mu}^{k}(E)$.    
The weighted $C^k$-spaces are Banach spaces but in general $C^{\infty}_{\mu}$ is not one.
\end{dfn}

We complete our collection of definitions by defining weighted
H\"older spaces. 

\begin{dfn}\label{s3dfn3}
Let $d(x,y)$ be the geodesic distance between points $x,y\in N$.  
Let $a\in (0,1)$, let $k\in\N$, let $\mu\in\R$ and let
\begin{align*}
H=\{(x,y)\in N&\times N\,:\,x\neq
y,\,c_1\rho(x)\leq\rho(y)\leq c_2\rho(x)\,\;\text{and}\;\,\\
&\text{there exists a geodesic in $N$ of length $d(x,y)$ from $x$ to
$y$}\},
\end{align*} 
where $0<c_1<1<c_2$ are constant. A section $\xi$
of $E$ is \emph{H\"older continuous} (with
\emph{exponent $a$}) if
$$[\xi]^a=\sup_{(x,y)\in H}\frac{|\xi(x)-\xi(y)|_{E}}{d(x,y)^a}<\infty.$$ 
We understand the numerator in the fraction above as follows.  Given
$(x,y)\in H$, there exists a 
geodesic $\gamma$ of length $d(x,y)$ connecting $x$ and $y$.
Parallel translation along $\gamma$ using the connection on $E$
identifies the fibres over $x$ and $y$ and the metrics on them.
Thus, with this identification, $|\xi(x)-\xi(y)|_E$ is well-defined.

The \emph{weighted H\"older space}
$C_{\mu}^{k,\,a}(E)$ is the
subspace of $C^{k,\,a}_{\text{loc}}(E)$ such
that the norm 
$$\|\xi\|_{C^{k,\,a}_{\mu}}=\|\xi\|_{C^{k}_{\mu}}+[\xi]^{k,\,a}_{\mu}$$
is finite, where $\|\,.\,\|_{C^k_{\mu}}$ is given in Definition \ref{s3dfn2} and
$$[\xi]^{k,\,a}_{\mu}=[\rho^{k+a-\mu}\nabla^k\xi]^{a}.
$$  Then $C_{\mu}^{k,\,a}(E)$ is a
Banach space.  It is clear that we have an embedding
$C_{\mu}^{k,\,a}(E)\hookrightarrow
C_{\mu}^l(E)$ whenever $l\leq k$.
\end{dfn}

\noindent The set $H$ in the definition above is introduced so that $[\xi]^{k,\,a}_{\mu}$ is well-defined.

\begin{note}
Given a weighted Banach space with weight $\mu$ on $N$, 
the dual weighted Banach space has weight $-3-\mu$.
\end{note}

We shall need the analogue of the Sobolev Embedding Theorem for
weighted spaces, which is adapted from \cite[Theorem 1.2]{Bartnik} and \cite[Lemma
7.2]{LockhartMcOwen}.
\begin{thm}[Weighted Sobolev Embedding Theorem]\label{s3thm4}
Let
 $p,q>1$,\\ $k,l\in\N$, $a\in (0,1)$ and $\mu,\nu\in\R$.
\begin{itemize}
\item[{\rm (i)}] If\/ $k\geq l$, $k-\frac{3}{p}\geq l-\frac{3}{q}$
and either
$p\leq q$ and $\mu\leq\nu$, or 
$p>q$ and $\mu<\nu$,
there is a continuous embedding
$L_{k,\,\mu}^p(E)\hookrightarrow$
$L_{l,\,\nu}^q(E)$. 
\item[{\rm (ii)}] If\/
$k-\frac{3}{p}\geq l+a$, there is a continuous embedding
$L_{k,\,\mu}^p(E)\hookrightarrow
C_{\mu}^{l,\,a}(E)$.
\end{itemize}
\end{thm}

We shall also require an Implicit Function Theorem for Banach
spaces, which follows immediately from \cite[Theorem 2.1]{Lang2}.
\begin{thm}[Implicit Function Theorem]\label{s3thm5}
$\!\!$Let $\mathcal{X}$ and $\mathcal{Y}$ 
be Banach spaces and let $\mathcal{U}\subseteq
\mathcal{X}$ be an open neighbourhood of $0$.  Let $\mathcal{F}:\mathcal{U}\rightarrow \mathcal{Y}$ be a
$C^k$-map $(k\geq 1)$ such that $\mathcal{F}(0)=0$.  Suppose further that
$d\mathcal{F}|_0:\mathcal{X}\rightarrow \mathcal{Y}$ is surjective with kernel $\mathcal{K}$ such
that $\mathcal{X}=\mathcal{K}\oplus \mathcal{Z}$ for some closed subspace $\mathcal{Z}$ of $\mathcal{X}$.  There
exist open sets $\mathcal{V}\subseteq \mathcal{K}$ and $\mathcal{W}\subseteq \mathcal{Z}$, both containing
$0$, with $\mathcal{V}\times \mathcal{W}\subseteq \mathcal{U}$, and a unique $C^k$-map
$\mathcal{G}:\mathcal{V}\rightarrow \mathcal{W}$ such that
$$\mathcal{F}^{-1}(0)\cap(\mathcal{V}\times \mathcal{W})=\big\{\big(x,\mathcal{G}(x)\big)\,:\,x\in \mathcal{V}\big\}$$
in $\mathcal{X}=\mathcal{K}\oplus \mathcal{Z}$.
\end{thm}

\begin{dfn}\label{subsetdfn}
For an open subset $U$ of $E$, $k\in\N$ and $\mu\in\R$, we introduce the notation $C^{k}_{\mu}(U)=\{\xi\in
C^{k}_{\mu}(E)\,:\,\Gamma_\xi\subseteq U\}$, where $\Gamma_{\xi}$ is the graph of $\xi$.  
We use similar conventions to define subsets of the weighted Sobolev and H\"older spaces, though we must take care that the spaces
 contain continuous sections, which then have well-defined graphs. 
\end{dfn}

\section[The deformation problem]{The deformation problem}\label{s4}

In this section we describe the deformation problem we are interested in 
and define a map whose kernel gives a local description of the moduli space.

\subsection[Setting up the problem]{Setting up the problem}\label{s4subs1}

We begin with the notation which we shall use throughout the rest of the paper. 
 We take $N$ to be an associative 3-fold in $\R^7$ 
which is asymptotically conical with rate $\lambda$ to a cone $C$ that has link $\Sigma$.  We use the notation of 
Definition \ref{s2dfn8} and recall that $C$ is associative by Proposition \ref{s2prop10}. 

\medskip

We define the moduli space for our deformation problem; that is, the space of AC associative 3-folds `near' $N$ in $\R^7$.

\begin{dfn}\label{s4dfn1}
The \emph{moduli space of deformations} $\mathcal{M}(N,\lambda)$ is
the set of associative 3-folds $N^\prime\subseteq \R^7$ which are
AC to $C$ with rate $\lambda$, and diffeomorphic to $N$, such that the inclusions $i:N\rightarrow\R^7$ and 
$i^{\prime}:N\cong N^{\prime}\rightarrow\R^7$ are isotopic through immersions $i^{\prime\prime}:N\rightarrow\R^7$,
where $i^{\prime\prime}(N)$ is AC to $C$ with rate $\lambda$ (but not necessarily associative).
\end{dfn}

\begin{note}
If $N$ is AC with rate $\lambda$ then it is also AC with rate $\lambda^{\prime}$ for all 
$\lambda^{\prime}\in[\lambda,1)$.  Thus, given any AC associative 3-fold $N$ we may define a one-parameter family of
moduli spaces of deformations of $N$.
\end{note}

The first result we need is immediate from the proof of
\cite[Chapter IV, Theorem 9]{Lang}.
\begin{thm}\label{s4thm2}
Let $P$ be a closed submanifold
of a Riemannian manifold $M$.  There exist an open subset $V$ of the
normal bundle $\nu(P)$ of $P$ in $M$, containing the zero section,
and an open set $S$ in $M$ containing $P$, such that the exponential
map $\exp|_V:V\rightarrow S$ is a diffeomorphism.
\end{thm}

\begin{note} The proof of this result relies entirely on the
observation that $\exp|_{\nu(P)}$ is a local isomorphism upon the
zero section.
\end{note} 

\noindent This information allows us to prove a corollary in a similar manner to \cite[Corollary 4.3]{Lotay1}, so we omit the proof.

\begin{cor}\label{s4cor3}
  Let
$C_{\infty}=\iota\big((R,\infty)\times\Sigma\big)$ and define
$n_\Psi: \nu(C_{\infty})\rightarrow\R^7$ by $$n_\Psi(r\sigma,v)=
v+\Psi(r,\sigma).$$ 
There exist an open subset $V$ of $\nu(C_{\infty})$, containing the
zero section, and an open set $S$ in $\R^7$ containing $N_{\infty}=N\setminus K$, such
that $n_\Psi|_{V}:V\rightarrow S$ is a diffeomorphism. 
Moreover, $V$ can be chosen such that $C^1_1(V)$ is an open neighbourhood of the zero section in $C^1_1\big(\nu(C_{\infty})\big)$.
\end{cor}

We now make an elementary observation which follows from our choice of diffeomorphism 
$\Psi:(R,\infty)\times\Sigma\rightarrow N_{\infty}$. 

\begin{lem}\label{s4lem4}
Use the notation of Corollary \ref{s4cor3} and of Definition \ref{s2dfn6} applied to $N$ and $C$.  Let $\mathcal{T}(C_{\infty})$ and $\mathcal{T}(N_{\infty})$ denote the obvious restrictions of the twisted spinor bundles over $C$ and $N$ respectively.  
There exist isomorphisms
$\upsilon:\nu(C_{\infty})\rightarrow\nu(N_{\infty})$ and $\Upsilon:\mathcal{T}(C_{\infty})\rightarrow\mathcal{T}(N_{\infty})$  
 such that the following diagram commutes:
\begin{equation*}
\begin{gathered}
\xymatrix{
\mathcal{T}(C_{\infty})\ar[d]^{\Upsilon}\ar[rr]^{\jmath_C^{-1}} &  &\nu(C_{\infty})\ar[d]^{\upsilon}\ar[rr] & &
C_{\infty}\ar[d]^{\Psi\circ\iota^{-1}}
\\
\mathcal{T}(N_{\infty})
\ar[rr]^{\jmath_N^{-1}} & &\nu(N_{\infty})\ar[rr] & &N_{\infty}.}
\end{gathered}
\end{equation*}
\end{lem}

\begin{proof}
By \eq{s2eq4}, we have a decomposition:
$$\R^7=T_{\Psi(r,\sigma)}N\oplus\nu_{r\sigma}(C)$$
whenever $r>R$.  Thus, $\nu_{\Psi(r,\sigma)}(N)\cong\nu_{r\sigma}(C)$ for all $r>R$, 
so $\nu(C_{\infty})$ and $\nu(N_\infty)$ are certainly isomorphic.  Moreover, again from \eq{s2eq4}, 
we have a distinguished section $\psi$ of the normal bundle $\nu(C_{\infty})$ essentially given
by the map $\Psi-\iota$ on $(R,\infty)\times\Sigma$.  Thus, if the projection from $\nu(N_\infty)$ to $N_\infty$ is 
$\pi$, the map $\psi\circ\iota\circ\Psi^{-1}\circ\pi$ is the inverse of the isomorphism $\upsilon$.  
We then let $\Upsilon=\jmath_N\circ\upsilon\circ\jmath_C^{-1}$.
\end{proof}

\begin{remark}
The fact that an isomorphism exists between $\nu(N_{\infty})$ and $\nu(C_{\infty})$, and thus between the associated twisted spinor bundles, 
is trivial.  The point of Lemma \ref{s4lem4} is that it enables us to relate the decay rate at infinity of sections of the 
isomorphic bundles.  This will be useful later.
\end{remark}

We can now conclude this subsection with our description of deformations of $N$ as the
graphs of `small' normal vector fields.

\begin{prop}\label{s4prop5} Use the notation of Corollary \ref{s4cor3} and Lemma \ref{s4lem4}.  
There exist an open subset $U$ of $\nu(N)$, containing the zero section and $\upsilon(V)$, 
a tubular neighbourhood $T$ of $N$ in $\R^7$, containing $S$, and a diffeomorphism $\delta:U\rightarrow T$,
affine on the fibres, such that $\delta(0)=N$ and the following diagram commutes:
\begin{equation}\label{s4eq1}
\begin{gathered}
\xymatrix{
\upsilon(V)\ar[rr]^{\upsilon^{-1}}\ar[dr]_{\delta} & &
V \ar[dl]^{n_\Psi}
\\
&S.&
}\end{gathered}
\end{equation}
Moreover, $U$ can be chosen such that $C^1_1(U)$ is an open neighbourhood of the zero section in $C^1_1\big(\nu(N)\big)$.
\end{prop}

\begin{proof}
We can define $\delta|_{\upsilon(V)}$ by \eq{s4eq1}, then smoothly extend $\upsilon(V)$, $S$ and $\delta|_{\upsilon(V)}$ over the 
compact set $K$ to $U$, $T$ and $\delta$ respectively, such that the conditions of the proposition are met.  Further, since
we can choose $C^1_1(V)$ to be open by Corollary \ref{s4cor3} and we
 have defined $\upsilon$ using $\Psi-\iota$, we can ensure that $U$ has the same property as $V$.
\end{proof}

\begin{remark}
The reader may wonder why we did not simply apply Theorem \ref{s4thm2} to all of $N$; that is, just consider exponential
deformations.  The reason is that Theorem \ref{s4thm2} does not gives us control of the asymptotic behaviour of our open sets and diffeomorphism, 
which we require for our analysis later to be valid.
\end{remark}

\subsection[The deformation map]{The deformation map}\label{s4subs2}

Recall the set $U$ given in Proposition \ref{s4prop5} and the notation in Definition \ref{subsetdfn}. 
The fact that $C^1_1(U)$ is open in $C^1_1(\nu(N))$ ensures that
$C^{k}_{\lambda}(U)$ is an open subset of
$C^{k}_{\lambda}(\nu(N))$ for $k\geq 1$, since $\lambda<1$. 
Moreover, the subsets of the weighted Sobolev and H\"older spaces will be open whenever the
space embeds continuously in $C^1_1$; for example, $L^p_{k+1,\,\lambda}(U)$ is an open neighbourhood of zero 
in $L^p_{k+1,\,\lambda}(\nu(N))$ if $k>\frac{3}{p}$ by Theorem \ref{s3thm4}(ii).  

Given the results of $\S$\ref{s4subs1}, we can define our \emph{deformation map}.

\begin{dfn}\label{s4dfn6}  Recall the form $\chi$ given in Definition \ref{s2dfn3}
and the notation of Definition \ref{subsetdfn} and Proposition \ref{s4prop5}. 
For $v\in C^1_{\text{loc}}( U)$, let
$g_{v}:N\rightarrow\Gamma_{v}\subseteq U$ be
given by $g_{v}(x)=\big(x,v(x)\big)$. Let
$f_v=\delta\circ g_{v}$ and let
$N_{v}=f_v(N)$, so that
$N_{v}\subseteq T$ is the deformation of $N$ corresponding to $v$.  
Since $N_v$ is essentially the graph of a normal vector field on $N$, the diffeomorphism $f_v:N\rightarrow N_v$ induces an isomorphism 
$f_v^*:\Omega^3(N_v;\nu(N_v))\rightarrow\Omega^3(N;\nu(N))$, abusing notation slightly. 
Then $f_v^*\big(\chi|_{N_{v}}\big)\in
C^{0}_{\text{loc}}\big(\Omega^3(N;\nu(N))\big)$ by Proposition \ref{s2prop5}, so we can
 define a nonlinear first-order differential operator 
$F:C^{1}_{\text{loc}}(U)\rightarrow
C^{0}_{\text{loc}}(\nu(N))$ by
\begin{equation*}\label{ch7s1eq1} 
F(v)=*f_v^*\big(\chi|_{N_{v}}\big). 
\end{equation*}
By Proposition \ref{s2prop4}, $\text{Ker}\,F$ is the set of $v\in C^1_{\text{loc}}(U)$ such that $N_v$ is associative.  

Following the ideas in the proof of \cite[Theorem 5.2]{McLean}, and using the notation of Definition \ref{normaldiffdfn} 
applied to $N\subseteq\R^7$, 
\begin{equation}\label{NDslasheq}
\Dslash v:=dF|_0(v)= *d^{\perp}_N(v\cdot\chi)|_N
\end{equation}
for $v\in C^1_{\text{loc}}(\nu(N))$.  In \eq{NDslasheq}, $(v\cdot\chi)|_N\in C^1_{\text{loc}}\big(\Omega^2(N;\nu(N))\big)$ by
 Proposition \ref{s2prop5}, so $d^{\perp}_N(v\cdot\chi)|_N\in C^0_{\text{loc}}\big(\Omega^3(N;\nu(N))\big)$ by Definition \ref
{normaldiffdfn} and hence $*d^{\perp}_N(v\cdot\chi)|_N\in C^0_{\text{loc}}(\nu(N))$.  Thus, $\Dslash$ is a well-defined first-order 
linear differential operator from $\nu(N)$ to itself.    

Recall the isomorphism $\jmath_N$ between $\nu(N)$ and a twisted spinor bundle 
$\mathcal{T}(N)$ on $N$ introduced
in Definition \ref{s2dfn6}.  In \cite[$\S$5]{McLean}, McLean computes the \emph{symbol} of the differential operator
 $\jmath_N\circ\Dslash\circ\jmath_N^{-1}:
\mathcal{T}(N)\rightarrow\mathcal{T}(N)$ and thus shows that it is a (twisted) 
\emph{Dirac operator}.  Therefore, $\Dslash=dF|_0$ is an \emph{elliptic} operator which can be thought of as a Dirac operator.
\end{dfn}

\begin{notes}
\begin{itemize}\item[]
\item[(a)] McLean \cite[$\S$5]{McLean} in effect only computes the symbol of $\Dslash$ because he does not take into account the fact that the exterior derivative does not act naturally on normal vector-valued forms on $N$.  One must choose a connection on $\nu(N)$ and thus $\Dslash$ contains a zero-order component which depends on this choice of
 connection.
\item[(b)] We can give an explicit local description of the first-order terms in $\Dslash$ as follows.  Recall Definition 
\ref{normaldiffdfn}
and the note after Definition \ref{s2dfn6}.  
Let $\nabla_i^{\perp}$ denote differentiation on $\nu(N)$ in the direction of the tangent vector $e_i$ on $N$, 
determined by $\nabla^{\perp}_N$.  Then, identifying $\nu(N)$ with $\mathcal{T}(N)$,    
the action of the first order terms in $\Dslash$ is locally given by: 
\begin{equation}\label{twistedDeq}
e_1\times\nabla_1^{\perp}+e_2\times\nabla_2^{\perp}+e_3\times\nabla_3^{\perp}.
\end{equation} 
\end{itemize}
\end{notes}

Though the kernel of $F$ in $C^1_{\text{loc}}(U)$ describes nearby associative deformations of $N$, we need to restrict to
deformations which are AC to $C$ with rate $\lambda$.

\begin{prop}\label{s4prop7}  Use the notation of Definitions \ref{s4dfn1} and \ref{s4dfn6}.  
The moduli space $\mathcal{M}(N,\lambda)$ is homeomorphic near $N$ to the kernel of $F:C^{\infty}_{\lambda}(U)\rightarrow 
C^{\infty}\big(\nu(N)\big)$.
\end{prop}

\begin{proof}  By Proposition \ref{s4prop5} and Definition \ref{s4dfn6}, every associative deformation of $N$ in the tubular 
neighbourhood $T$ can be written as $N_v$ for some $v\in C^1_{\text{loc}}(U)$ such that $F(v)=0$.  Moreover, since deformations $N_v$  
must necessarily be smooth, $v\in C^{\infty}(U)$.  Recall the notation from Definition \ref{s2dfn8}. 
Now, $N_v$ is AC with rate $\lambda$ to $C$ 
if and only if there exist a compact subset $K_v$ of $N_v$ and a diffeomorphism $\Psi_v:(R,\infty)\times\Sigma
\rightarrow N_v\setminus K_v$ satisfying \eq{s2eq4} and \eq{s2eq5}.  By condition \eq{s2eq4}, $\Psi-\iota$ and $\Psi_v-\iota$ 
define normal vector fields on $C_{\infty}=\iota((R,\infty)\times\Sigma)$ and hence, by Lemma \ref{s4lem4}, elements $w$ and $w_v$ in 
$\nu(N_{\infty})$ respectively.  Moreover, by \eq{s2eq5},
$w$ and $w_v$ lie in $C^{\infty}_{\lambda}(\nu(N_{\infty}))$.  Finally, we notice that $w_v=v+w$ on $N_{\infty}$, which implies 
that $v\in C^{\infty}_{\lambda}(U)$ as required.
\end{proof}

It follows from our next proposition that $F$ defines a smooth map between weighted Banach spaces, which is key for the
 proof of Theorem \ref{thm1}.  
The result can be derived using the method presented in the proof of \cite[Proposition 4.9]{Lotay1}, so we refer the
 interested reader to that article.

\begin{prop}\label{s4prop8} Apply Definition \ref{normaldiffdfn} to $N$.
The map $F$, given in Definition \ref{s4dfn6}, can be written as
\begin{equation}\label{s4eq2}
F(v)(x)=\Dslash v(x)+P_F\big(x,v(x),\nabla^{\perp}_N v(x)\big)
\end{equation}
for $x\in N$ and $v\in C^1_{\text{\emph{loc}}}(U)$, where 
$$P_F:\{(x,y,z)\,:\,(x,y)\in U,\,z\in
T_x^*N\otimes\nu_x(N)\}\rightarrow 
\nu(N)$$ 
is a
smooth map such that $P_F(x,y,z)\in 
\nu_x(N)$.  Denote $P_F\big(x,v(x),\nabla^{\perp}_N v(x)\big)$ by
$P_F(v)(x)$ for $x\in N$.  

Let $\mu<1$ and suppose $v\in C^1_1(U)$ with $\|v\|_{C^1_1}$ sufficiently small.  
For each $n\in\N$, if\/ $v \in C^{n+1}_\mu(U)$, $P_F(v)\in
C^n_{2\mu-2}\big(\nu(N)\big)$ and there exists a constant $c_n>0$
such that
\begin{equation*}
\|P_F(v)\|_{C_{2\mu-2}^n}\leq
c_n\|v\|_{C^{n+1}_{\mu}}^2\,.
\end{equation*}
Let $k\geq 2$, $l\in\N$ and $a\in(0,1)$.  
If\/ $v\in L^2_{k+1,\,\mu}(U)$ or $v\in C^{l+1,\,a}_{\mu}(U)$ with $\|v\|_{C^1_1}$ sufficiently small, 
$P_F(v)\in L^2_{k,\,2\mu-2}\big(\nu(N)\big)$ or 
$C^{l,\,a}_{2\mu-2}\big(\nu(N)\big)$ and there exist constants $l_{k}>0$ and $c_{l,\,a}>0$ such that 
$$\|P_F(v)\|_{L^2_{k,\,2\mu-2}}\leq l_{k}\|v\|_{L^2_{k+1,\,\mu}}^2
\quad\text{or}\quad \|P_F(v)\|_{C^{l,\,a}_{2\mu-2}}\leq c_{l,\,a}\|v\|_{C^{l+1,\,a}_{\mu}}^2.$$
\end{prop}

\begin{remark}
As $\mu<1$, $2\mu-2<\mu-1$, so $C^n_{2\mu-2}\hookrightarrow C^n_{\mu-1}$.  Similar  
continuous embeddings occur for the weighted Sobolev and H\"older spaces.  Furthermore, the conditions 
$k\geq 2$ and $\mu<1$ ensure, by Definition \ref{s3dfn3} and Theorem \ref{s3thm4}(ii), that 
$C^{l+1,\,a}_{\mu}\hookrightarrow C^1_1$ and $L^2_{k+1,\,\mu}\hookrightarrow C^1_1$.
\end{remark}

Proposition \ref{s4prop8} helps us to prove an important regularity result for elements of the kernel of $F$.  

\begin{prop}\label{s4prop9}  Use the notation of Definition \ref{s4dfn6}.
Let $v\in L^2_{k+1,\,\mu}(U)$ for some $k\geq 3$ and $\mu<1$ satisfy $F(v)=0$.  If $\|v\|_{C^1_1}$ is sufficiently small, then 
$v\in C^{\infty}_{\mu}(U)$.
\end{prop}

\begin{proof}
Suppose $v\in C^{l-1,\,\frac{1}{2}}_{\mu}(U)$ for some $l\geq 3$.  
Let $G(v)=\Dslash F(v)$ and, for convenience, set $\nabla=\nabla^{\perp}_N$ as given by Definition \ref{normaldiffdfn} 
applied to $N$.  Since $G(v)$ is linear in $\nabla^2v$ with coefficients depending on $v$ and $\nabla v$, we can write
$$G(v)(x)=R_F\big(x,v(x),\nabla v(x)\big)\nabla^2v(x)+E_F\big(x,v(x),\nabla v(x)\big),$$
where $R_F$ and $E_F$ are smooth functions of their arguments.  We use the notation $R_F(v)$ and $E_F(v)$ in a similar way to
$P_F(v)$ in Proposition \ref{s4prop8}.  Having done this, we can define a \emph{linear} second-order differential operator on 
$C^2_{\text{loc}}(\nu(N))$, with coefficients depending on $v$ and $\nabla v$, by
$$S_v(w)=R_F(v)\nabla^2w.$$
Notice that $S_v$ is \emph{not} the linearisation of $G$ at $0$, which is $\Dslash^2$.  However, 
$S_v$ defines an \emph{elliptic} operator $S_v:C^{l-1,\,\frac{1}{2}}_{\mu}(\nu(N))\rightarrow C^{l-3,\,\frac{1}{2}}_{\mu-2}(\nu(N))$
with coefficients in $C^{l-2,\,\frac{1}{2}}_{\text{loc}}$.

Suppose further that $F(v)=0$ and $\|v\|_{C^1_1}$ is sufficiently small. 
Then $G(v)=0$, so $S_v(v)=-E_F(v)$.  However, $E_F(v)$ depends only on $v$ and $\nabla v$ and is at
worst quadratic in these quantities by Proposition \ref{s4prop8}, so $E_F(v)\in C^{l-2,\,\frac{1}{2}}_{\mu-2}(\nu(N))$ as we 
have control of the first $l-1$ derivatives of $v$.  Thus, by the Schauder elliptic regularity estimates given in \cite{Mazya}, 
$v\in C^{l,\,\frac{1}{2}}_{\mu}(U)$.  

Therefore, if $v\in L^2_{k+1,\,\mu}(U)$ with $F(v)=0$, $v\in C^{k-1,\,\frac{1}{2}}_{\mu}(U)$ by 
Theorem \ref{s3thm4}(ii), so we can use induction
to prove that $v\in C^{l-1,\,\frac{1}{2}}_{\mu}(U)$ for all $l\geq k$. 
\end{proof}

\section[The Dirac operator]{The Dirac operator}\label{s5}

In this section we discuss the Fredholm and index theory of the operator $\Dslash$
 using the scheme of \cite{LockhartMcOwen}. 
We begin with our first Fredholm theory result, which follows from \cite[Theorem 1.1 \& Theorem 6.1]{LockhartMcOwen}.

\begin{prop}\label{s5prop1} There exists a countable discrete set $\mathcal{D}\subseteq\R$ such that $\Dslash$, given by
Definition \ref{s4dfn6}, defines a Fredholm operator $\Dslash: L^2_{k+1,\,\mu}\big(\nu(N)\big)\rightarrow 
L^2_{k,\,\mu-1}\big(\nu(N)\big)$, for any
$k\in\N$, if and only if $\mu\notin\mathcal{D}$.  Moreover, $\mathcal{D}$ is determined by a first-order differential operator 
acting on $\Sigma$.
\end{prop}

\noindent As in \cite[$\S$5]{Lotay1} and \cite[$\S$6.1.2]{Marshall}, we can describe this set $\mathcal{D}$ explicitly.  

We begin by defining the Dirac operator on $C$.

\begin{dfn}\label{CDiracdfn} 
By applying Definition \ref{normaldiffdfn} to the cone $C$, we can define a first-order differential operator $\Dslash_C$ on 
$\nu(C)$ by $\Dslash_Cv_C=*d^{\perp}_C(v_C\cdot\chi)|_C$.  This is well-defined since $C$ is associative by Proposition \ref{s2prop10}.  
\end{dfn} 

\begin{note}
The operator $\Dslash$ given by Definition \ref{s4dfn6} is asymptotically conical (with rate $\lambda$)
 to $\Dslash_C$ in the sense described in \cite[$\S$4.3]{Lotay1}.
\end{note} 

We now observe that if $\tau$ is a tensor of type $(p,q)$ on the cone $C$ which is independent of the radial coordinate $r$,  then, using index notation, 
$$|\tau|^2=\tau^{a_1,\ldots,a_p}_{c_1,\ldots,c_q}\tau^{b_1,\ldots,b_p}_{d_1,\ldots,d_q}g_{a_1b_1}\ldots
g_{a_pb_p}g^{c_1d_1}\ldots g^{c_qd_q},$$
where $g=g_\text{cone}$ is the cone metric on $C$.  Thus, since the cone metric is of order $O(r^2)$, we see that $\tau$ is of order $O(r^{p-q})$ on $C$.  

This calculation allows us to deduce that a normal vector field $v_\Sigma$ on the link $\Sigma$, pulled back to $C$, is naturally 
of order $O(r)$.  Hence, $v_C(r\sigma)=r^{\mu-1}v_\Sigma(\sigma)$ is a homogeneous normal vector field on $C$ of 
order $O(r^{\mu})$.  

Since the operator $\Dslash$ on $N$ is asymptotic to the operator $\Dslash_C$ on the cone $C$, the work in 
\cite{LockhartMcOwen} implies that $\mu\in\mathcal{D}$ given by Proposition \ref{s5prop1} if and only if there exists a
 nonzero homogeneous solution $v_C$ to $\Dslash_Cv_C=0$ of order $O(r^{\mu})$.  Thus, we are compelled to study 
the equation 
\begin{equation}\label{dperpeq1}
d_C^{\perp}(r^{\mu-1}v_\Sigma\cdot\chi)|_C=0.
\end{equation}
Applying Definition \ref{normaldiffdfn} to $\Sigma\subseteq\mathcal{S}^6$, we get a derivative $d_\Sigma^{\perp}$ 
on normal vector-valued forms on $\Sigma$ and we notice that 
\begin{equation}\label{dperpeq2}
d_C^{\perp}=d_\Sigma^{\perp}+dr\w\frac{\partial}{\partial r}.
\end{equation}
Using the fact that $C$ is associative, we can perform a straightforward local calculation using the formulae \eq{s2eq3} and \eq{chiconeeq} for $\chi$ to deduce that 
\begin{equation}\label{dperpeq3}
(r^{\mu-1}v_\Sigma\cdot\chi)|_C= (r^{\mu+1}Jv_\Sigma)\vol_\Sigma-r^{\mu}dr\wedge (v_\Sigma\cdot\zeta)|_\Sigma,
\end{equation}
where $J$ is the almost complex structure on $\mathcal{S}^6$ given in Definition \ref{6spheredfn} and 
$(v_\Sigma\cdot\zeta)|_\Sigma$ is a normal vector-valued 1-form on $\Sigma$.  Notice that since $\Sigma$ is 
a pseudoholomorphic curve by Lemma \ref{phololem}, $J$ preserves the normal spaces of $\Sigma$ and so $Jv_\Sigma$ is 
a normal vector field on $\Sigma$.
Substituting \eq{dperpeq2}-\eq{dperpeq3} into \eq{dperpeq1} gives:
\begin{equation}\label{dperpeq4}
d_C^{\perp}(r^{\mu-1}v_\Sigma\cdot\chi)|_C=0\quad \Leftrightarrow\quad *_\Sigma d_\Sigma^{\perp}(v_\Sigma\cdot\zeta)|_\Sigma=-(\mu+1)Jv_\Sigma,
\end{equation}
where $*_\Sigma$ is the Hodge star on $\Sigma$.
These calculations encourage us to make the following definition.

\begin{dfn}\label{SigmaDiracdfn}
Apply Definition \ref{normaldiffdfn} to $\Sigma\subseteq\mathcal{S}^6$ and recall the notation of 
Definition \ref{6spheredfn} and the form $\zeta$ given by \eq{chiconeeq}.  Define an operator $\Dslash_\Sigma$
on $v_{\Sigma}\in C^{\infty}\big(\nu(\Sigma)\big)$ via 
$\Dslash_\Sigma v_\Sigma= *_\Sigma d_\Sigma^{\perp}(v_\Sigma\cdot\zeta)|_\Sigma$.  This is a well-defined first-order linear differential 
operator from $\nu(\Sigma)$ to itself.  The isomorphism between $\nu(C)$ and a twisted spinor bundle given in Definition 
\ref{s2dfn6} induces an isomorphism between $\nu(\Sigma)$ and 
a twisted spinor bundle $\mathcal{T}(\Sigma)$ over $\Sigma$.  Thus, $\Dslash_\Sigma$ can be viewed as a 
Dirac operator on $\mathcal{T}(\Sigma)\cong\nu(\Sigma)$.  
\end{dfn}

\begin{note}
As in the note after Definition \ref{s2dfn6}, at $r\sigma\in C$ we can identify $T_{r\sigma}\R^7$ with $\Im\O$ such 
that $T_{r\sigma}C=\langle e_1,e_2,e_3\rangle$.  Moreover, since $\GG_2$ acts transitively on the set of 2-planes 
in $\R^7$, we can further adapt frames so that $e_3$ is the radial vector and $\langle e_1,e_2\rangle\cong T_\sigma\Sigma$.
Notice that taking cross product with $e_3$ defines the almost complex structure on $\Sigma$ by 
Definition \ref{6spheredfn}.  Using \eq{twistedDeq} and \eq{dperpeq1}-\eq{dperpeq4}, we deduce that the first-order 
terms in $\Dslash_\Sigma$ have the local formula:
\begin{equation}\label{twistedDeq2}
e_1\times\nabla_1^{\perp}+e_2\times\nabla_2^{\perp},
\end{equation}
where $\nabla_i^{\perp}$ is now determined by the connection $\nabla_\Sigma^{\perp}$ on $\nu(\Sigma)$ induced from 
the Levi--Civita connection of the round metric on $\mathcal{S}^6$ as in Definition \ref{normaldiffdfn}.
\end{note}

We take a moment to make an important observation about $\Dslash_\Sigma$.  

\begin{lem}\label{antilem} In the notation of Definitions \ref{6spheredfn} and \ref{SigmaDiracdfn},  $\Dslash_\Sigma(Jv_\Sigma)=-J(\Dslash_\Sigma v_\Sigma)$ for all sections $v_\Sigma$ of $\nu(\Sigma)$.
\end{lem}

\begin{proof} 
By the note after Definition \ref{SigmaDiracdfn}, it is now a straightforward local calculation using multiplication on $\Im\O$
and \eq{twistedDeq2} to see that $J$ and the first-order terms in $\Dslash_\Sigma$ anti-commute, so we need only focus on 
the zero-order piece of $\Dslash_\Sigma$. 

Recall that $\Sigma$ is a pseudoholomorphic curve in $\mathcal{S}^6$ with respect to $J$.  
Let $\{f_1,f_2,f_3\}$ be a local unitary frame for $T^{1,0}\mathcal{S}^6|_{\Sigma}$ adapted such that $f_1$ spans $T^{1,0}\Sigma$, 
and let $\theta_1$ be the dual 1-form to $f_1$.  
Then, we can locally write a normal vector field $v_\Sigma$ on $\Sigma$ as $\Re(w_2f_2+w_3f_3)$
for suitable functions $w_2$ and $w_3$.  A straightforward computation leads us to deduce that, locally, 
\begin{equation}\label{cxeq1}
(v_\Sigma\cdot\zeta)|_\Sigma=\Im(w_3\bar{f}_2\theta_1-w_2\bar{f}_3\theta_1).
\end{equation}  
  
Since $\Sigma$ is a pseudoholomorphic curve in $\mathcal{S}^6$, it follows from the structure equations given in \cite[$\S$4]{Bryant} that there exist real 1-forms $\rho_0,\rho_1,\rho_2,\rho_3$ on $\Sigma$ such that:
\begin{align}
d^{\perp}_\Sigma f_2&=i f_2\rho_2+f_3(\rho_0+i\rho_1)-\bar{f}_3\theta_1;\label{cxeq2}\\
d^{\perp}_\Sigma f_3&=f_2(-\rho_0+i\rho_1)+if_3\rho_3+\bar{f}_2\theta_1.\label{cxeq3}
\end{align}
Calculation using \eq{cxeq1}-\eq{cxeq3} leads us to deduce that the component of $\Dslash_\Sigma v_\Sigma$ containing no derivatives 
of $w_2$ and $w_3$ is
 $\Re(u_2f_2+u_3f_3)$ where $u_2$ and $u_3$ are linear in $\bar{w}_2$ and $\bar{w}_3$.   Since the action of $J$ on $v_\Sigma$ is
 equivalent to the map $(w_2,w_3)\mapsto (iw_2,iw_3)$, we see that $J(\Dslash_\Sigma v_\Sigma)$ has the
 opposite sign to $\Dslash_\Sigma(Jv_\Sigma)$ as claimed.
\end{proof}

Lemma \ref{antilem} motivates the next definition.

\begin{dfn}\label{dbardfn}
Recall Definitions \ref{6spheredfn} and \ref{SigmaDiracdfn}.  Let 
$$\bar{\partial}_\Sigma=J\Dslash_\Sigma:C^{\infty}\big(\nu(\Sigma)\big)\rightarrow C^{\infty}\big(\nu(\Sigma)\big).$$    
Since $J$ acts on $\nu(\Sigma)$ as an almost complex structure, $\nu(\Sigma)$ can be viewed as a complex vector bundle.  The 
operator $\bar{\partial}_{\Sigma}$ anti-commutes with $J$ by Lemma \ref{antilem} and so can be thought of as a $\bar{\partial}$-operator 
on $\nu(\Sigma)$, which then defines a holomorphic structure on $\nu(\Sigma)$. 
\end{dfn}

\begin{note}
Following on from the note after Definition \ref{SigmaDiracdfn}, we see that the local formula for the first-order terms in $\bar{\partial}_\Sigma$ is:
\begin{equation*}\label{twistedDeq3}
e_1\times\nabla_2^{\perp}-e_2\times\nabla_1^{\perp}.
\end{equation*}
\end{note}

Returning to the Fredholm theory of $\Dslash$, we deduce from Definitions \ref{CDiracdfn} and \ref{SigmaDiracdfn} and
\eq{dperpeq1}-\eq{dperpeq4} that  
\begin{equation}\label{Diraceq0}
\Dslash_C v_C=\frac{1}{r^3}\left(
\Dslash_\Sigma+Jr\frac{\partial}{\partial r}\right)(r^2v_C).
\end{equation}
Therefore, $\Dslash_C v_C=0$ for $v_C=r^{\mu-1}v_\Sigma$ if and only if 
\begin{equation}\label{Diraceq1}
\Dslash_{\Sigma}v_{\Sigma}=-(\mu+1)Jv_{\Sigma}.
\end{equation}
By Definition \ref{dbardfn}, \eq{Diraceq1} is equivalent to
\begin{equation}\label{Diraceq2}
\bar{\partial}_\Sigma v_\Sigma=(\mu+1)v_\Sigma.
\end{equation}
By the theory in \cite{LockhartMcOwen} we can derive our description of $\mathcal{D}$.

\begin{prop}\label{s5prop2}
For $\mu\in\R$, let
$$D(\mu)=\big\{v_{\Sigma}\in C^{\infty}\big(\nu(\Sigma)\big)\,:\,\bar{\partial}_{\Sigma}v_{\Sigma}=(\mu+1)v_{\Sigma}\big\},$$
where $\bar{\partial}_{\Sigma}$ is described in Definition \ref{dbardfn}. 
 In the notation
of Proposition \ref{s5prop1}, $\mathcal{D}=\{\mu\in\R\,:\,D(\mu)\neq 0\}$.
\end{prop}

We can now give our index result which follows from \cite[Theorem 1.2]{LockhartMcOwen}.

\begin{prop}\label{s5prop3}
Use the notation of Definition \ref{s4dfn6} and Proposition \ref{s5prop1}.  
Let $\d(\mu)$ be the dimension of the space of
solutions $v_C$ to $\Dslash_C v_C=0$ of the form
$$v_C(r\sigma)=r^{\mu-1}v_{\Sigma}(r,\sigma),$$
where $v_{\Sigma}(r,\sigma)$ is a polynomial in $\log r$ with coefficients in $C^{\infty}\big(\nu(\Sigma)\big)$.  
Let $\ind_{\mu}(\Dslash)$ denote the Fredholm index of $\,\Dslash: L^2_{k+1,\,\mu}\big(\nu(N)\big)\rightarrow
L^2_{k,\,\mu-1}\big(\nu(N)\big)$.  If $\lambda,\lambda^{\prime}\in\R\setminus\mathcal{D}$ such that $\lambda\leq\lambda^{\prime}$,
$$\ind_{\lambda^{\prime}}(\Dslash)-\ind_{\lambda}(\Dslash)=\sum_{\mu\in\mathcal{D}\,
\cap(\lambda,\,\lambda^{\prime})}\!\!\!\!\!\!\d(\mu).$$ 
\end{prop}

We give an alternative characterisation of the quantity $\d(\mu)$ by proving that
there can be no $\log r$ terms in solutions $v_C$ to $\Dslash_C v_C=0$ of the form given in Proposition \ref{s5prop3}.  
This is an analogous result to \cite[Proposition 5.6]{Lotay1}.

\begin{prop}\label{s5prop4} In the notation of Propositions \ref{s5prop2}-\ref{s5prop3},
$\d(\mu)=\dim D(\mu)$.
\end{prop}

\begin{proof}
Let $m\in\N$ and let
$$v_C(r\sigma)=r^{\mu-1}\sum_{j=0}^m (\log r)^jv_j(\sigma),$$
where $v_j\in C^{\infty}(\nu(\Sigma))$ and $v_m\neq 0$, be such that $\,\Dslash_Cv_C=0$.  Thinking of $\,\Dslash_Cv_C=0$ as a polynomial equation in
$\log r$ and comparing coefficients of $(\log r)^m$, we see from \eq{Diraceq0} that
$$\Dslash_{\Sigma}v_m=-(\mu+1)Jv_m.$$

Suppose, for a contradiction, that $m\geq 1$.  By comparing coefficients of $(\log r)^{m-1}$ in $\Dslash_Cv_C=0$ using 
\eq{Diraceq0} again, we have that
$$mJv_m+\Dslash_{\Sigma}v_{m-1}+(\mu+1)Jv_{m-1}=0.$$
Noticing that $\Dslash_\Sigma$ is formally self-adjoint and using Lemma \ref{antilem}, we calculate:
\begin{align*}
m\langle v_m,v_m\rangle_{L^2}&=
\langle J\Dslash_{\Sigma}v_{m-1},v_m\rangle_{L^2}-(\mu+1)\langle v_{m-1},v_m\rangle_{L^2}\\
&=\langle v_{m-1}, \Dslash_{\Sigma}(-Jv_m)-(\mu+1)v_m\rangle_{L^2}\\
&=\langle v_{m-1}, J\Dslash_{\Sigma}v_m-(\mu+1)v_m\rangle_{L^2}=0.
\end{align*}
Therefore, $v_m=0$, our required contradiction.

Hence, $v_C(r\sigma)=r^{\mu-1}v_0(\sigma)$ and $v_0\in D(\mu)$.  The result follows.
\end{proof}

Suppose $\mu\in\mathcal{D}$ and that $v_{\Sigma}$ is a normal vector field on $\Sigma$ satisfying \eq{Diraceq2}.  By Lemma \ref{antilem}, $\bar{\partial}_{\Sigma}(Jv_\Sigma)=-(\mu+1)Jv_\Sigma$ and so $Jv_\Sigma$ is also an eigenvector
of $\bar{\partial}_{\Sigma}$, but now with eigenvalue $-(\mu+1)=(-2-\mu)+1$.  In fact, the map $v\mapsto Jv$ defines an isomorphism between $D(\mu)$ and $D(-2-\mu)$ so that $\d(\mu)=\d(-2-\mu)$.  A standard result from the theory of Dirac operators on compact manifolds (c.f.~\cite[$\S$4.2]{Friedrich})
 states that the spectrum of such an operator consists of a sequence of eigenvalues $\lambda_n$ such that 
$|\lambda_n|\rightarrow\infty$.  We may infer that $\bar{\partial}_\Sigma$ 
 has infinitely many eigenvalues symmetrically situated above and below $-1$.

We collect these observations together in the following lemma.

\begin{lem}\label{s5lem1} Using the notation of Proposition \ref{s5prop3},   
 $$\mathcal{D}\setminus\{-1\}=\left\{-1\pm\mu_n\,:\,n\in\mathbb{N}\right\}$$ such that $\mu_n>0$ for all $n$
and $\mu_n\rightarrow\infty$ as $n\rightarrow\infty$.  Moreover, $\d(-1+\mu_n)=\d(-1-\mu_n)$ for all $n$.
\end{lem}

\begin{notes}
We can interpret some other distinguished points of $\mathcal{D}$ geometrically.
\begin{itemize}
\item[(a)]  The rate $0$ corresponds to the introduction of deformations of $N$
given by translation.  Thus, $0\in\mathcal{D}$ and, if $m$ is the number of ends of $N$ which are non-planar,
$$\d(0)\geq 7m+4\big(b^0(\Sigma)-m\big)=4b^0(\Sigma)+3m>0.$$
Moreover, if $n$ is the number of components of $N$ which are non-planar, translations of these components 
give genuine deformations of $N$, so at least $4b^0(N)+3n$ is added to the kernel of $\Dslash$ at $0$.  Explicitly, 
these normal vector fields are simply the projections onto $\nu(N)$ of the translation vector fields on $\R^7$.
\item[(b)]  If we were to allow $N$ to be asymptotically conical to $C$ with rate $1$ or higher, then $N$ would be 
AC to any deformation of $C$.  Therefore, rate $1$ corresponds to deformations of the link $\Sigma$ of $C$.  
Hence we can interpret solutions to 
$$\bar{\partial}_{\Sigma}v_{\Sigma}=2v_{\Sigma}$$
as infinitesimal deformations of $\Sigma$ as a pseudoholomorphic curve in $\mathcal{S}^6$.
\item[(c)] The normal vector fields added at rate $-1$ are in the kernel of $\bar{\partial}_\Sigma$ and so can be 
thought of as holomorphic sections of $\nu(\Sigma)$.  Equivalently, these vector fields lie in the kernel of $\Dslash_\Sigma$ and so 
can be viewed as harmonic twisted spinors.  Notice that they do \emph{not} define infinitesimal deformations of 
$\Sigma$ as a pseudoholomorphic curve in $\mathcal{S}^6$ as one might naively expect.
\end{itemize}
\end{notes}

\section[The moduli space]{The moduli space}\label{s6}

In this section we prove our main results.  We begin by identifying the infinitesimal deformation space and the obstruction space.

\begin{dfn}\label{s6dfn1}
Let $k\geq 3$ and recall the notation of Definition \ref{s4dfn6}.  The \emph{infinitesimal deformation space} is
$$\mathcal{I}(N,\lambda)=\big\{v\in L^2_{k+1,\,\lambda}\big(\nu(N)\big)\,:\,\Dslash v=0\big\}.$$
By elliptic regularity results as given in \cite{Mazya}, $\mathcal{I}(N,\lambda)$ is finite-dimensional and independent of $k$.
\end{dfn}

\begin{dfn}\label{s6dfn2}
Use the notation of Definition \ref{s4dfn6} and Proposition \ref{s5prop1}.
Let $k\geq 3$ and suppose that $\lambda\notin\mathcal{D}$.  
Since $$\Dslash:L^2_{k+1,\,\lambda}\big(\nu(N)\big)\rightarrow L^2_{k,\,\lambda-1}\big(\nu(N)\big)$$ is Fredholm, it has closed range and finite-dimensional cokernel.  Thus, there exists a finite-dimensional subspace
$\mathcal{O}(N,\lambda)$ of $L^2_{k,\,\lambda-1}(\nu(N))$ such that
$$L^2_{k,\,\lambda-1}(\nu(N))=\Dslash\big(L^2_{k+1,\,\lambda}(\nu(N))\big)\oplus\mathcal{O}(N,\lambda).$$
We call $\mathcal{O}(N,\lambda)$ the \emph{obstruction space}.

Moreover, $\mathcal{O}(N,\lambda)$ is isomorphic to the kernel of the adjoint map
$$\Dslash^*=\Dslash:L^2_{l+1,\,-2-\lambda}\big(\nu(N)\big)\rightarrow L^2_{l,\,-3-\lambda}\big(\nu(N)\big),$$ for any $l\in\N$, because the dual space 
of the Sobolev space with weight $\lambda-1$ has weight $-3-(\lambda-1)=-2-\lambda$,
$\Dslash$ is self-adjoint and its kernel is independent of $l$ by elliptic regularity.  Further, we 
can choose $\mathcal{O}(N,\lambda)$ to equal this kernel if $\lambda\geq -1$, since then $L^2_{k+1,\,-2-\lambda}\hookrightarrow 
L^2_{k+1,\,\lambda-1}$ by Theorem \ref{s3thm4}(i) as $-2-\lambda\leq\lambda-1$.
\end{dfn}

\begin{remark}
We should really consider the adjoint operator acting on a Sobolev space with $l=-k-1$, but 
because we are only interested in its kernel, which is independent of $l\in\Z$ by elliptic regularity, 
we can avoid this technicality.
\end{remark} 

Having made these definitions, we can prove the following theorem.

\begin{thm}\label{s6thm3}
Let $N$ be an associative 3-fold in $\R^7$ which is asymptotically conical with rate $\lambda$ to a cone $C$ with link $\Sigma$. 
Use the notation of Definitions \ref{s4dfn1}, \ref{s6dfn1} and \ref{s6dfn2} and Proposition \ref{s5prop1}. 
Suppose further that $\lambda\notin\mathcal{D}$.  

There exists a smooth manifold $\hat{\mathcal{M}}(N,\lambda)$, which is 
an open neighbourhood of $0$ in $\mathcal{I}(N,\lambda)$, and a smooth map 
$\pi:\hat{\mathcal{M}}(N,\lambda)\rightarrow\mathcal{O}(N,\lambda)$, with $\pi(0)=0$, such that $\mathcal{M}(N,\lambda)$ near $N$ 
is homeomorphic to $\pi^{-1}(0)$ near $0$.  Thus, the expected dimension of the moduli space is
 $\dim\mathcal{I}(N,\lambda)-\dim\mathcal{O}(N,\lambda)$ 
 and $\mathcal{M}(N,\lambda)$ is smooth if $\mathcal{O}(N,\lambda)=\{0\}$.
\end{thm}

\begin{remark}
More generally, if $\pi$ is the zero map then $\mathcal{M}(N,\lambda)$ 
is smooth with dimension equal to that of $\mathcal{I}(N,\lambda)$.
\end{remark}

\begin{proof}
For some $k\geq 3$, let 
\begin{align*}
\mathcal{X}=L^2_{k+1,\,\lambda}\big(\nu(N)\big)\times\mathcal{O}(N,\lambda)\quad
\text{and}\quad
\mathcal{Y}= L^2_{k,\,\lambda-1}\big(\nu(N)\big).
\end{align*}
Then $\mathcal{X}$ and $\mathcal{Y}$ are clearly Banach spaces as $\mathcal{O}(N,\lambda)$ is a finite-dimensional subspace of a 
Banach space.  Let $$\mathcal{U}= 
L^2_{k+1,\,\lambda}(U)\times\mathcal{O}(N,\lambda),$$ where $U$ is given by 
Proposition \ref{s4prop5}.  As noted at the start of $\S$\ref{s4subs2}, $\mathcal{U}$ is an open neighbourhood of zero in $\mathcal{X}$.  Moreover,  by making $U$ smaller,
 which we are free to do, we can decrease the $C^1_1$-norm of all $v\in C^1_1(U)$ so that the relevant 
requirements of Propositions \ref{s4prop8} and \ref{s4prop9} are met by $v$.

We now define a smooth map of 
Banach spaces $\mathcal{F}:\mathcal{X}\rightarrow\mathcal{Y}$ by 
$$\mathcal{F}(v,w)=F(v)+w,$$
where $F$ is given in Definition \ref{s4dfn6}.  This map is well-defined by Proposition \ref{s4prop8} and the choice of $U$ described
above. 
Since $N$ is associative, $F(0)=0$, so $\mathcal{F}(0,0)=0$.  Moreover, 
$d\mathcal{F}|_{(0,0)}:\mathcal{X}\rightarrow\mathcal{Y}$
 acts as
$$d\mathcal{F}|_{(0,0)}:(v,w)\mapsto \Dslash v +w.$$
By the definition of the obstruction space in Definition \ref{s6dfn2}, $d\mathcal{F}|_{(0,0)}$ is surjective and 
$d\mathcal{F}|_{(0,0)}(v,w)=0$ if and only if $(\Dslash v,w)=(0,0)$.  Thus, the kernel of $d\mathcal{F}|_{(0,0)}$, 
$\mathcal{K}=
\mathcal{I}(N,\lambda)\times\{0\}$, is finite-dimensional and splits $\mathcal{X}$; that is, there exists a closed
subspace $\mathcal{Z}$ of $\mathcal{X}$ such that $\mathcal{K}\oplus\mathcal{Z}=\mathcal{X}$.  Moreover, we can write
$\mathcal{Z}=
\mathcal{Z}_1\times\mathcal{O}(N,\lambda)$ for a closed subspace $\mathcal{Z}_1$ of 
$L^2_{k+1,\,\lambda}(\nu(N))$.

Applying Theorem \ref{s3thm5}, we have open sets $\mathcal{V}\subseteq\mathcal{I}(N,\lambda)$, 
$\mathcal{W}_1\subseteq\mathcal{Z}_1$ and $\mathcal{W}_2\subseteq\mathcal{O}(N,\lambda)$, and 
smooth maps $\mathcal{G}_j:\mathcal{V}\rightarrow\mathcal{W}_j$ for $j=1,2$ such that
$$\mathcal{F}^{-1}(0)\cap \big((\mathcal{V}\times\mathcal{W}_1)\times\mathcal{W}_2\big)=
\big\{\big((v,\mathcal{G}_1(v)),\mathcal{G}_2(v)\big)\,:\,v\in\mathcal{V}\big\}$$
in $\mathcal{X}=(\mathcal{I}(N,\lambda)\oplus\mathcal{Z}_1)\times\mathcal{O}(N,\lambda)$.  That is, the kernel of $\mathcal{F}$ is 
diffeomorphic near $(0,0)$ to an open subset of $\mathcal{I}(N,\lambda)$ containing $0$.  

Now define $\hat{\mathcal{M}}(N,\lambda)=\mathcal{V}$ and $\pi:\hat{\mathcal{M}}(N,\lambda)\rightarrow \mathcal{O}(N,\lambda)$ by
 $\pi(v)=\mathcal{G}_2(v)$.  Thus, 
an open neighbourhood of zero in $F^{-1}(0)$ is homeomorphic to an open neighbourhood of zero in $\pi^{-1}(0)$.  Finally, 
by Proposition \ref{s4prop9} (which is applicable by our choice of $U$), $F^{-1}(0)$ consists of smooth forms as $k\geq 3$.  The result follows from Proposition \ref{s4prop7}.
\end{proof}
\begin{remarks}
\begin{itemize}\item[]
\item[(a)] When $\lambda\in\mathcal{D}$ or $\mathcal{O}(N,\lambda)\neq\{0\}$ the moduli space may not be smooth, 
or may have larger than expected dimension.
\item[(b)] $\mathcal{I}(N,\lambda)\neq \{0\}$ whenever $N$ is non-planar, since dilations of $N$ are examples of 
AC associative deformations of $N$.  This is proved in Proposition \ref{s6prop8}.
\item[(c)] Recall from Lemma \ref{s5lem1} that $\mathcal{D}$ has infinitely many elements above and below $-1$.  
Thus, $\mathcal{O}(N,\lambda)\neq\{0\}$ for $\lambda$ sufficiently negative and zero for $\lambda$ sufficiently 
above $-1$.
\end{itemize}
\end{remarks}

Theorem \ref{s6thm3} proves Theorem \ref{thm1}. To prove Theorem \ref{thm2}, we calculate 
the expected dimension of the moduli space using the index theory
in $\S$\ref{s5}. We start with some notation for convenience.

\begin{dfn}\label{s6dfn4}
Let $k,l\in\N$, let $\mu\in\R$ and use the notation of Definition \ref{s4dfn6}.  Let
$$\Dslash_{\mu}=\Dslash: L^2_{k+1,\,\mu}\big(\nu(N)\big)\rightarrow L^2_{k,\,\mu-1}\big(\nu(N)\big).$$
The formal adjoint of this operator is
$$\Dslash_{\mu}^*=\Dslash: L^2_{l+1,\,-2-\mu}\big(\nu(N)\big)\rightarrow L^2_{l,\,-3-\mu}\big(\nu(N)\big),$$
as noted in Definition \ref{s6dfn2}. 
 Let $\mathcal{K}(\mu)$ and $\mathcal{C}(\mu)$ be the kernels of $\Dslash_\mu$ and
$\Dslash^*_\mu$.  Notice that $\mathcal{K}(\lambda)=\mathcal{I}(N,\lambda)$ and, if $\lambda\notin\mathcal{D}$,
 $\mathcal{C}(\lambda)\cong\mathcal{O}(N,\lambda)$.  Moreover,  
we can choose $\mathcal{O}(N,\lambda)=\mathcal{C}(\lambda)$ if $\lambda\geq -1$ and $\lambda\notin\mathcal{D}$.  
 In this notation, the expected dimension of $\mathcal{M}(N,\lambda)$ is the index of $\Dslash_{\lambda}$, denoted 
 $\ind_{\lambda}(\Dslash)$.
\end{dfn}

From this definition we have an elementary lemma.

\begin{lem}\label{s6lem5}
Use the notation of Propositions \ref{s5prop1} and \ref{s5prop3} and Definition \ref{s6dfn4}.  
Let $\epsilon>0$ be such that
$[-1-\epsilon,-1+\epsilon]\cap\mathcal{D}\subseteq\{-1\}$.
\begin{itemize}
\item[\emph{(i)}] 
 $\mathcal{K}(-1)=\mathcal{C}(-1)$ and hence, if $-1\notin\mathcal{D}$, $\ind_{-1}(\Dslash)=0$.
\item[\emph{(ii)}] If $-1\in\mathcal{D}$, 
 $\ind_{-1-\epsilon}(\Dslash)=-\frac{\d(-1)}{2}$ and $\ind_{-1+\epsilon}
(\Dslash)=\frac{\d(-1)}{2}$.  
In particular, $\ind_\mu(\Dslash)>0$ for $\mu>-1$.
\item[\emph{(iii)}] $\ind_{\mu}(\Dslash)<0$ for $\mu<-2$.
\end{itemize}
\end{lem}

\begin{remark}
Notice that the formulae in (ii) imply $\d(-1)$ is even.  This tallies with our interpretation of $\d(-1)
=\dim\Ker\bar{\partial}_\Sigma$ as the (real) dimension of the space of holomorphic sections of $\nu(\Sigma)$. 
\end{remark}

\begin{proof}
Part (i) immediately follows from Definition \ref{s6dfn4} and the observation that $-2-(-1)=-1$.  

For part (ii), first note that
$$\ind_{-1+\epsilon}(\Dslash)\geq \dim\mathcal{K}(-1)-\dim\mathcal{C}(-1)
\geq \ind_{-1-\epsilon}(\Dslash),$$ because of the inclusions $L^2_{k+1,\,-1-\epsilon}\hookrightarrow
L^2_{k+1,\,-1}\hookrightarrow L^2_{k+1,\,-1+\epsilon}$.  We also know that, by Proposition \ref{s5prop3}, 
$$\ind_{-1+\epsilon}(\Dslash)-\ind_{-1-\epsilon}(\Dslash)=\d(-1)>0.$$  Now we note that 
$$\dim\mathcal{K}(-1\pm\epsilon)=\dim\mathcal{C}(-1\mp\epsilon),$$ simply by comparing the 
operators in Definition \ref{s6dfn4}.  It follows that 
\begin{align*}
\dim\mathcal{K}(-1+\epsilon)&=\dim\mathcal{K}(-1-\epsilon)+\frac{\d(-1)}{2}\quad\text{and}\\
\dim\mathcal{C}(-1+\epsilon)&=\dim\mathcal{C}(-1-\epsilon)-\frac{\d(-1)}{2}.
\end{align*}
The result is an easy consequence of these formulae.

For (iii), by Lemma \ref{s5lem1} and the remarks following it we know that $-2\in\mathcal{D}$, 
and therefore $\ind_{\mu}(\Dslash)\leq -\d(-2)<0$ for all $\mu<-2$, irrespective of whether $-1$ is in
$\mathcal{D}$ or not.
\end{proof}

From this result we can quickly deduce our formulae for the expected dimension of the moduli space simply by using 
Proposition \ref{s5prop3}.  This completes the proof of Theorem \ref{thm2}.

\begin{prop}\label{s6prop6}
Use the notation of Propositions \ref{s5prop1} and \ref{s5prop3}.  Let $\mu_{-}$ be the greatest element of $\mathcal{D}$ which is strictly less than $-1$. 
\begin{itemize}
\item[\emph{(i)}]
If $\lambda\in(-1,1)\setminus\mathcal{D}$,
$$\ind_{\lambda}(\Dslash)=
\frac{1}{2}\,\d(-1)+\!\!\!\!\sum_{\mu\in\mathcal{D}\,
\cap(-1,\,\lambda)}\!\!\!\!\!\!\d(\mu).$$
In particular, the index is positive for $\lambda>0$.
\item[\emph{(ii)}] 
If $\lambda\in (-\infty,-1]\setminus\mathcal{D}$, $\ind_{\lambda}(\Dslash)\leq -\frac{1}{2}\,\d(-1)$, and  
if $\lambda\in (-\infty,\mu_{-})\setminus\mathcal{D}$, $\ind_{\lambda}(\Dslash)<0$.  In particular, 
the index is negative for $\lambda<-2$.
\end{itemize}
\end{prop}

We conclude that the virtual dimension of the moduli space is non-negative if $\lambda>-1$,
 $\lambda\notin\mathcal{D}$.  
If $\bar{\partial}_{\Sigma}$, or equivalently $\Dslash_\Sigma$, 
 acting on the normal bundle of the link $\Sigma$ of the asymptotic cone $C$ has a non-zero kernel, 
then $\mathcal{M}(N,\lambda)$ should be a positive-dimensional manifold for generic rates $\lambda>-1$.  We can think of this 
 condition as requiring that $\nu(\Sigma)$ admits holomorphic sections in the sense of Definition \ref{dbardfn}.

Notice that the dimension of the moduli space is determined entirely by the link $\Sigma$ so that, 
unlike in the coassociative and special Lagrangian cases, there is no topological component from the
associative 3-fold in the calculation of the expected dimension.  
This tallies with the fact that compact associative 3-folds are expected to be isolated, whereas 
compact coassociative and special Lagrangian submanifolds occur in finite-dimensional families based on their
topology.  

\begin{note}
One could mirror all of these results, with minor changes, for associative 3-folds in $\R^7$ which are 
\emph{asymptotically cylindrical}.  However, one would discover in this case that the index of the Dirac operator, 
for all relevant growth rates, is non-positive.  
Therefore we would expect, generically, that asymptotically cylindrical 
associative 3-folds are isolated, as in the compact case.
\end{note}

We now provide some bounds on the dimensions of the infinitesimal deformation and obstruction spaces.

\begin{prop}\label{s6prop7}
Use the notation of Definitions \ref{s4dfn6}, \ref{s6dfn1} and \ref{s6dfn2} and Proposition \ref{s5prop1}.  Let 
$$\mathcal{H}(N)= \big\{v\in L^2\big(\nu(N)\big)\,:\,\Dslash v=0\big\}.$$  
\begin{itemize}
\item[\emph{(i)}] If $\lambda\geq -\frac{3}{2}\geq\lambda^{\prime}$,
$\dim\mathcal{I}(N,\lambda)\geq \dim\mathcal{H}(N)\geq\dim
\mathcal{I}(N,\lambda^{\prime})$.
\item[\emph{(ii)}] If $\mu\geq -\frac{1}{2}\geq\mu^{\prime}$ and $\mu,\mu^{\prime}\notin\mathcal{D}$, 
$\dim\mathcal{O}(N,\mu)\leq \dim\mathcal{H}(N)\leq\dim\mathcal{O}(N,\mu^{\prime})$.  
\end{itemize}
\end{prop}

\begin{proof}  Recall that $L^2=L^2_{0,\,-\frac{3}{2}}$ and that solutions to $\Dslash v=0$ in any weighted
Sobolev space are smooth. 
We deduce (i) from Definition \ref{s6dfn1} and the inclusions $L^2_{0,\,\lambda^{\prime}}\hookrightarrow 
L^2_{0,\,-\frac{3}{2}}\hookrightarrow L^2_{0,\,\lambda}$
given by Theorem \ref{s3thm4}(i).  
Part (ii) follows similarly from Definition \ref{s6dfn2} and the 
inclusions $L^2_{0,\,-2-\mu}\hookrightarrow L^2_{0,\,-\frac{3}{2}}\hookrightarrow 
L^2_{0,\,-2-\mu^{\prime}}$.
\end{proof}

We now prove an interesting result which uses the fact that there is always a deformation 
of an AC associative 3-fold given by dilation.

\begin{prop}\label{s6prop8}
Use the notation of Definition \ref{s6dfn1} and 
Propositions \ref{s5prop1} and \ref{s6prop7}.  Let $n$ be the number of components of $N$ which are non-planar.  
\begin{itemize}
\item[\emph{(i)}] For any $\lambda<1$, $\lambda\notin\mathcal{D}$, 
$\dim\mathcal{I}(N,\lambda)\geq n$.  If, in addition, 
$\lambda>0$, $\dim\mathcal{I}(N,\lambda)\geq 4(b^0(N)+n)$.
\item[\emph{(ii)}] Suppose $\mathcal{H}(N)=0$.  If $\lambda<-\frac{3}{2}$, $n=0$; that is, $N$ is a union of 3-planes.  
\end{itemize}
\end{prop}

\begin{proof}  For convenience, suppose that $N$ is connected so that $n=0$ or $n=1$.

We can always define an element of $C^{\infty}_{\lambda}(\nu(N))$ corresponding to an infinitesimal deformation of $N$ as follows.   Let $u=\sum_{i=1}^7x_i\frac{\partial}{\partial x_i}$ 
be the dilation vector field on $\R^7$ and let $v_u$ be the orthogonal projection of $u$ onto $\nu(N)$, so that $tv_u$ is the
normal vector field associated to the deformation $N\mapsto e^tN$ for $t\in\R$.  Clearly, $e^tN$ is also AC with rate $\lambda$ to $C$
and 
$v_u\in C^{\infty}_{\lambda}(\nu(N))$.  
Moreover, using the formula \eq{s2eq3} for $\chi$, the fact that $\chi$ vanishes on 
$TN$ and Definition \ref{s4dfn6},
$$\Dslash v_u=*d^{\perp}_N(v_u\cdot\chi)|_N=*d^{\perp}_N(u\cdot\chi)|_N=3*\!\chi|_N=0.$$

For (i), first let $\epsilon>0$ be such that $[\lambda,\lambda+\epsilon]\cap\mathcal{D}=\emptyset$, which is
possible as $\mathcal{D}$ is discrete and $\lambda\notin\mathcal{D}$.  By Proposition \ref{s5prop3}, there are 
no changes in $\ind_{\mu}(\Dslash)$ for $\mu\in[\lambda,\lambda+\epsilon]$, hence no changes in the kernel or
 cokernel of $\Dslash_{\mu}$ for these rates $\mu$.  Since $C^{\infty}_{\lambda}\hookrightarrow L^2_{k+1,\,\lambda+\epsilon}$ 
for any $k\in\N$ and  $\mathcal{I}(N,\lambda+\epsilon)=\mathcal{I}(N,\lambda)$, $v_u\in\mathcal{I}(N,\lambda)$.  
If $n=1$, so that $N\neq\R^3$, $v_u\neq 0$ and therefore defines a 1-dimensional subspace of $\mathcal{I}(N,\lambda)$.  
The second sentence in (i) follows from note (a) after Lemma \ref{s5lem1}, which states that
at least $4+3n$ is added to $\mathcal{I}(N,\lambda)$ as the rate crosses $0$ and the additional normal 
vector fields are given by translations.  

For (ii), we know that if $\lambda<-\frac{3}{2}$, $C^{\infty}_{\lambda}\hookrightarrow L^2_{0,\,-\frac{3}{2}}=L^2$, 
so that $v_u\in\mathcal{H}(N)$.  However, $\mathcal{H}(N)=0$, so $v_u=0$.  
Thus, $N$ is a cone in $\R^7$, so it is a linear $\R^3$ as it is nonsingular.

These arguments can be easily adapted to the case where $N$ is not connected by considering each component separately.
\end{proof}

We conclude with a simple, but interesting, corollary which follows from Theorem \ref{s6thm3} and Propositions \ref{s6prop7} and 
\ref{s6prop8}. 

\begin{cor}\label{s6cor9}
Use the notation of Theorem \ref{s6thm3} and Proposition \ref{s6prop7}.  If $\lambda\geq-\frac{1}{2}$, 
$\mathcal{H}(N)=0$ 
and $N$ has a non-planar component, then 
$N$ admits a smooth moduli space of deformations $\mathcal{M}(N,\lambda)$ with positive dimension.
\end{cor}

\section[Examples]{Examples}\label{s7}

In this final section, we begin by presenting examples of asymptotically conical associative 3-folds in $\R^7$ 
which do not arise from lower-dimensional geometries.  We also consider the cases where AC associative 3-folds 
are formed from calibrated submanifolds in $\C^3$. 

\subsection[Ruled associative 3-folds]{Ruled associative 3-folds}

In \cite[$\S$4]{Lotay2}, the author produced examples of associative 3-folds in $\R^7$ which are asymptotically conical 
to cones invariant under an action of $\U(1)$.  Moreover, these examples are \emph{ruled};  
that is, fibred by affine straight lines over a surface $\Sigma$.  It was noted in \cite[$\S$6]{Lotay3} that a ruled associative 3-fold is AC with rate $0$ or lower 
to a naturally associated cone with link $\Sigma$, if $\Sigma$ is compact.

The author studied ruled associative 3-folds in detail in \cite[$\S$6]{Lotay3} and gave methods for constructing them 
starting from associative cones.  In particular, it is observed that the underlying surface $\Sigma$ 
admits a natural complex
 structure, and that one may use a holomorphic vector field on $\Sigma$ to construct a ruled 
 associative 3-fold diffeomorphic to $\Sigma\times\R$.
 For this technique to be implemented, one must start with a associative cone over $\mathcal{S}^2$ or $T^2$.  

\medskip
 
 The author 
\cite[$\S$4]{Lotay2} 
was able to produce examples of associative $T^2$-cones using elementary techniques, from which ruled associative 
3-folds could then be constructed.  

\begin{remark}
There are many examples of associative $T^2$-cones given by special 
Lagrangian torus cones in $\C^3$.  The holomorphic vector field technique applied to these cones will only produce
 (ruled) AC special Lagrangian 3-folds.  
\end{remark}

We begin by describing a family of $\U(1)$-invariant associative cones in a result which follows 
from \cite[Theorems 4.2-4.4]{Lotay2}.

\begin{thm}\label{s7thm1} 
Let $x_1:\R\rightarrow\R$ and $z_1,z_2,z_3:\R\rightarrow\C$ be smooth functions such that:
\begin{gather*}
x_1^2+|z_1|^2+|z_2|^2+|z_3|^2=1;\qquad \Im(z_1)=0;\\
\Re(z_1z_2z_3)=a_1;\qquad
|z_1|(x_1^2+|z_1|^2-1)=a_2; 
\\
\Re\!\big(z_1(z_2^2-z_3^2)\big)=a_3;\quad\text{and}\quad
\Im\!\big(z_1(z_2^2+z_3^2)\big)=a_4
\end{gather*}
for some real constants $a_j$.  For generic choices of $a_j$, 
$$C=\big\{\big(rx_1(t),\,re^{2is}z_1(t),\,re^{-is}z_2(t) ,
\,re^{-is}z_3(t)\big):\,r>0,\,s,t\in\R\big\}$$
is an associative cone in $\R^7\cong\R\oplus\C^3$, which is invariant under the $\U(1)$ action
$(\xi_1,\zeta_1,\zeta_2,\zeta_3)\mapsto (\xi_1,e^{2is}\zeta_1,e^{-is}\zeta_2,e^{-is}\zeta_3)$, and has link 
diffeomorphic to $T^2$.
\end{thm}

\noindent We can now present examples of AC associative 3-folds \cite[Theorem 4.7]{Lotay2} which 
are formed by applying the holomorphic vector field construction to the $T^2$-cones of Theorem \ref{s7thm1}.

\begin{thm}\label{s7thm2}  Use the notation of Theorem \ref{s7thm1}.
Let $u,v:\R^2\rightarrow\R$ be
functions satisfying the Cauchy--Riemann equations.  The subset $N_{(u,v)}$ of\/ $\R\oplus\C^3$
given by
\begin{align}N_{(u,v)}=
\bigg\{\bigg(&rx_1(t)+v(s,t)\Big(2|z_1(t)|^2-|z_2(t)|^2-|z_3(t)|^2\Big),\nonumber\\
&\,e^{2is}\Big(r+2iu(s,t)-2v(s,t)x_1(t)\Big)z_1(t),\nonumber\\&\,e^{-is}\Big(\big(r-iu(s,t)+v(s,t)x_1(t)
\big)z_2(t)
-3iv(s,t)\overline{z_3(t)z_1(t)}\,\Big),\nonumber\\
&\,e^{-is}\Big(\big(r-iu(s,t)+v(s,t)x_1(t)\big)z_3(t)+3iv(s,t)\overline{z_1(t)z_2(t)}\,\Big)
\bigg)\,:\,\nonumber\\
&\qquad\qquad\qquad\qquad\qquad\qquad\qquad\qquad\qquad\quad
\,r>0,\,s,t\in\R\bigg\}\label{Nuv}\end{align} is an
associative 3-fold in $\R^7\cong\R\oplus\C^3$ which
is asymptotically conical to $C$ with rate $-1$.
\end{thm}

\begin{remark}
$N_{(u,v)}$ above will be nonsingular if $u,v$ are not identically zero, but will not, in general, be $\U(1)$-invariant 
even though $C$ is.  
It would be interesting to know if there is a critical rate below which an AC associative 3-fold inherits 
the symmetries of its asymptotic cone, as occurs in the coassociative and special Lagrangian scenarios.
\end{remark}

First, as $-1$ is a possible critical decay rate, we must consider $\lambda>-1$.  If 
$\lambda$ is such that $(-1,\lambda]\cap\mathcal{D}=\emptyset$, the expected dimension of the moduli space
$\mathcal{M}(N_{(u,v)},\lambda)$ is $\frac{1}{2}\,\d(-1)\geq 0$ from Theorem \ref{s6thm3} and Proposition \ref{s6prop6}.  
We see that the dilation $N_{(u,v)}\mapsto e^\tau N_{(u,v)}$ sends $N_{(u,v)}$ to $N_{(u,e^\tau v)}$, where the
latter is defined by \eq{Nuv} for the pair of functions $(u,e^\tau v)$.  Clearly $(u,e^\tau v)$ satisfy the Cauchy--Riemann
equations if and only if $\tau=0$ (unless $u$ and $v$ are constant).  Thus, for $u,v$ non-constant, 
the dilations of $N_{(u,v)}$ do not lie in the family given by Theorem \ref{s7thm2} and are AC associative deformations
with rate $-1$.  If $u,v$ are instead both non-zero constants, then $N_{(u,v)}\cong N_{(u,1)}$ via dilation.  
In both cases we deduce that $\d(-1)\geq 2$.

If $\lambda>0$ and $\lambda\notin\mathcal{D}$, then the dimension of $\mathcal{M}(N_{(u,v)},\lambda)$ is expected to be at least $8$ 
by Proposition \ref{s6prop8}(i).
 
\begin{note} 
Bryant \cite[$\S$4]{Bryant} showed that there are many associative cones over $\mathcal{S}^2$, in contrast to 
the special Lagrangian case where there are no nontrivial $\mathcal{S}^2$-cones \cite[Theorem B]{Haskins}. A particular example is  
 given by the cone over the Bor\r{u}vka sphere $\mathcal{S}^2(\frac{1}{6})$ in $\mathcal{S}^6$, described in \cite{Boruvka}.  Given any 
 associative cone over $\mathcal{S}^2$, the holomorphic vector field construction will produce a family of 
AC associative 3-folds depending on the $6$-dimensional space of holomorphic vector fields on $\mathcal{S}^2$.
\end{note}

\subsection[Special Lagrangian 3-folds]{Special Lagrangian 3-folds}

We noted in Lemma \ref{s2lem1} that special Lagrangian 3-folds in $\C^3$ are generalised by associative 3-folds 
in $\R^7$.  Thus, by embedding $\C^3$ as a hyperplane in $\R^7$ we can produce examples of AC associative 
3-folds from special Lagrangian ones.  
However, the natural question is: are there any 
AC associative deformations of an embedded AC special Lagrangian 3-fold $L$ in $\R^7$ 
which are \emph{not} embedded special Lagrangian 3-folds?  Here, we give a partial answer to this question.

We begin with a definition.

\begin{dfn}\label{ACSLdfn}
Let $\R^7=\R\oplus\C^3$ with $x_1$ the coordinate on $\R$, and $\omega$ and $\Omega$ the 
K\"ahler and holomorphic volume form on $\C^3$ respectively, as in Definition \ref{SLdfn}.  Let $L$ be a special Lagrangian 3-fold in $\C^3\subseteq\R^7$.  Suppose further that $L$ is AC with rate 
$\lambda$ to a special Lagrangian cone $C$ with (minimal Legendrian) link $\Sigma$ in 
$\mathcal{S}^5\subseteq\mathcal{S}^6$.  Let $\mathcal{M}(L,\lambda)$ and $\tilde{\mathcal{M}}(L,\lambda)$ be 
the associative and special Lagrangian moduli spaces of AC deformations of $L$ respectively.

If $\Delta_\Sigma$ is the Laplacian on functions on $\Sigma$, we also define 
$$\tilde{\mathcal{D}}=\{\mu\in\R\,:\,\text{$(\mu+1)(\mu+2)$ is an eigenvalue of $\Delta_\Sigma$}\}$$
and let $\tilde{\d}(\mu)$ be the multiplicity of the eigenvalue $(\mu+1)(\mu+2)$ for $\mu\in\tilde{\mathcal{D}}$.
\end{dfn}

A step in the direction of answering our question is the following.  Suppose $L$ is AC with rate $\lambda<0$. 
The moduli space $\mathcal{M}(L,\lambda)$ consists of AC associative 3-folds $L^{\prime}$
 which are
also AC with rate $\lambda<0$ to $C$.  By \cite[Theorem 5.10]{Lotay3}, such $L^{\prime}$ must 
themselves be (embedded) special Lagrangian.  Therefore, by the deformation theory of AC special Lagrangian 3-folds 
\cite[Theorem 8.3.10]{Joyce3}, we deduce the following.

\begin{prop}\label{lowratesprop}  Use the notation of Definition \ref{ACSLdfn}.  If $\lambda<0$ then 
  $\tilde{\mathcal{M}}(L,\lambda)=\mathcal{M}(L,\lambda)$.  In particular, 
\begin{itemize}
\item[\emph{(i)}] if $\lambda\in(-1,0)\setminus\tilde{\mathcal{D}}$, then 
$\mathcal{M}(L,\lambda)$ is a smooth manifold of dimension
$$b^1(L)-b^0(L)+\!\!\!\!\sum_{\mu\in(-1,\lambda)\cap{\tilde{\mathcal{D}}}}\!\!\!\!\tilde{\d}(\mu);$$
\item[\emph{(ii)}] if $\lambda\in(-2,-1)$, then $\mathcal{M}(L,\lambda)$ is a smooth manifold of dimension $b^2(L)$.
\end{itemize}
\end{prop} 

\begin{remarks}  
Proposition \ref{lowratesprop}(ii) may seem to be at odds with Proposition \ref{s6prop6}(ii), where we showed that for
 $\lambda\in(-2,-1)$ the expected dimension of $\mathcal{M}(L,\lambda)$ was non-positive.  There is no contradiction, however, 
 since $L$ is special Lagrangian and so certainly not a generic example of an AC associative 3-fold.  Moreover, the condition 
 $\lambda\notin\tilde{\mathcal{D}}$ is not needed in Proposition \ref{lowratesprop}(ii) since $\tilde{\mathcal{D}}\cap(-2,-1)=\emptyset$.  
\end{remarks}

\noindent From Proposition \ref{lowratesprop}, it is clear that we must consider rates $\lambda>0$ if we are to get any truly associative deformations of $L$.  

It is obvious that the infinitesimal associative deformation space for 
$\lambda>0$ is automatically larger than the special Lagrangian one because we can translate in the $x_1$
 direction.  Thus we should consider the quotient of $\mathcal{M}(L,\lambda)$ by this translation action.  
However, the author has been unable to show whether there are further 
infinitesimal associative deformations for these rates which extend to genuine associative deformations.  Therefore, 
although it is possible that $\mathcal{M}(L,\lambda)/\R$ is strictly larger than $\tilde{\mathcal{M}}(L,\lambda)$, the problem remains unresolved.

\medskip

Proposition \ref{lowratesprop} has rate $-2$ as its strict lower bound for discussing deformation 
theory.  However, there are a number of non-planar special Lagrangian 3-folds in $\C^3$ which are AC with rate $-2$,  
described in \cite[Example 8.3.15]{Joyce3}.  
 We notice from Proposition \ref{s6prop8}(ii) 
that the space of $L^2$ solutions to $\Dslash v=0$ must be non-zero for these AC special Lagrangian 3-folds.

\subsection[Holomorphic curves]{Holomorphic curves}

Inspired by Lemma \ref{s2lem1}, we could consider submanifolds of $\R^7=\R\oplus\C^3$ of the form $N=\R\times S$, 
for some holomorphic curve $S$ in $\C^3$.  For $N$ to be an AC associative 3-fold we would naturally require $S$ to be
 asymptotically conical to a cone in $\C^3$.  
However, $N$ would then be asymptotically conical to a cone which
did not have an isolated singularity at $0$.  Thus, our deformation theory would not apply in this case.  

\medskip
 
\begin{funding}
National Science Foundation Mathematical Sciences Postdoctoral Research Fellowship.
\end{funding} 
 
\begin{ack}
The author would like to thank Dominic Joyce for helpful criticisms and suggestions, and MSRI for hospitality during 
the time of this project.  The author also thanks the referee for providing detailed and useful comments and suggestions. 
\end{ack}





\begin{thebibliography}{99}

\bibitem{Bartnik} R.~Bartnik, {\it The Mass of an Asymptotically
Flat Manifold}, Comm.~Pure Appl.~Math.~{\bf 39} (1986), 661--693.

\bibitem{Boruvka} O.~Bor\r{u}vka, {\it Sur les Surfaces Represent\'ees par les Fonctions Sph\'eriques de Premiere Esp\'ece}, J.~Math.~Pure et Appl.~{\bf 12} (1933), 337--383. 

\bibitem{Bryant} R.~L.~Bryant, {\it Submanifolds and Special Structures on the Octonions}, J.~Differential Geom.~{\bf 17}
 (1982), 185--232.

\bibitem{Friedrich} T.~Friedrich, {\it Dirac Operators in Riemannian Geometry}, Graduate Studies in Mathematics {\bf 25}, 
Amer. Math. Soc., Providence, Rhode Island, 2000.

\bibitem{Gayet} D.~Gayet and F.~Witt, {\it Deformations of Associative Submanifolds with Boundary}, preprint, arXiv:0802.1283.

\bibitem{HarLaw} R.~Harvey and H.~B.~Lawson, {\it Calibrated Geometries},
Acta Math.~{\bf 148} (1982), 47--152.

\bibitem{Haskins} M.~Haskins, {\it Special Lagrangian Cones}, Amer.~J.~Math.~{\bf 126} (2004), 845--871.  

\bibitem{Joyce5} D.~D.~Joyce, {\it Special Lagrangian Submanifolds
with Isolated Conical Singularities. I. Regularity}, Ann.~Global
Anal.~Geom.~{\bf 25} (2004), 201--251.  

\bibitem{Joyce6} D.~D.~Joyce, {\it Special Lagrangian Submanifolds
with Isolated Conical Singularities. II. Moduli Spaces}, Ann.~Global
Anal.~Geom.~{\bf 25} (2004), 301--352.   

\bibitem{Joyce7} D.~D.~Joyce, {\it Special Lagrangian Submanifolds
with Isolated Conical Singularities. III. Desingularization, The
Unobstructed Case}, Ann.~Global Anal.~Geom.~{\bf 26} (2004), 1--58.

\bibitem{Joyce8} D.~ D.~Joyce, {\it Special Lagrangian Submanifolds
with Isolated Conical Singularities. IV. Desingularization,
Obstructions and Families}, Ann.~Global Anal.~Geom.~{\bf 26} (2004),
117--174.  

\bibitem{Joyce9} D.~D.~Joyce, {\it Special Lagrangian Submanifolds
with Isolated Conical Singularities. V. Survey and Applications}, J.~Differential 
Geom.~{\bf 63} (2003), 299--347.   

\bibitem{Joyce3} D.~D.~Joyce, {\it Riemannian Holonomy Groups and Calibrated Geometry}, 
Oxford Graduate Texts in Mathematics {\bf 12}, OUP, Oxford, 2007.  

\bibitem{JoySalur} D.~D.~Joyce and S.~Salur, {\it Deformations of Asymptotically Cylindrical Coassociative
Submanifolds with Fixed Boundary}, Geom.~Topol.~{\bf 9} (2005), 1115--1146.  

\bibitem{Lang} S.~Lang, {\it Differential Manifolds},
Addison-Wesley, Reading, Massachusetts, 1972.

\bibitem{Lang2} S.~Lang, {\it Real Analysis}, Addison-Wesley, Reading, Massachusetts, 1983.

\bibitem{LockhartMcOwen} R.~B.~Lockhart and R.~C.~McOwen, {\it
Elliptic Differential Operators on Noncompact Manifolds}, Ann.~Sc.~Norm.~Super.~Pisa Cl.~Sci.~{\bf 12} (1985), 409--447.

\bibitem{Lotay3} J.~Lotay, {\it Constructing Associative 3-folds by Evolution
Equations},\\ Comm.~Anal.~Geom.~{\bf 13} (2005), 999--1037.  

\bibitem{Lotay2} J.~D.~Lotay, \textit{Calibrated Submanifolds of $\R^7$ and $\R^8$ with Symmetries},\\ 
Q.~J.~Math.~{\bf 58} (2007), 53--70.  

\bibitem{Lotay1} J.~D.~Lotay, \textit{Deformation Theory of
Asymptotically Conical Coassociative 4-folds}, Proc.~London~Math.~Soc~{\bf 99} (2009), 386--424.

\bibitem{Marshall} S.~P.~Marshall, {\it Deformations of Special
Lagrangian Submanifolds},\\ D.Phil.~thesis, Oxford University, Oxford, 2002.  

\bibitem{Mazya} V.~G.~Maz'ya and B.~Plamenevskij, {\it Elliptic Boundary Value Problems},\\ Amer.~Math.~Soc.~Transl.~{\bf 123} (1984), 1--56.

\bibitem{McLean} R.~C.~McLean, {\it Deformations of Calibrated
Submanifolds},\\ Comm.~Anal.~Geom.~{\bf 6} (1998), 705--747.

\bibitem{Pacini} T.~Pacini, {\it Deformations of Asymptotically
Conical Special Lagrangian Submanifolds}, Pacific J.~Math.~{\bf
215} (2004), 151--181.  

\bibitem{Salur} S.~Salur, {\it Deformations of Asymptotically Cylindrical Coassociative Submanifolds with Moving Boundary},
preprint, arXiv:math/0601420.

\end{thebibliography}
\end{document}